\documentclass[reqno]{amsart}
\usepackage{graphicx,helvet}
\usepackage{epstopdf}
\usepackage[percent]{overpic}
\usepackage{xcolor}
\usepackage{multirow}

\usepackage{epsfig}
\usepackage{cite}
\usepackage[english]{babel}
\usepackage[utf8]{inputenc}

\allowdisplaybreaks

\renewcommand{\epsilon}{\varepsilon}

\newcommand{\R}{\mathbb{R}}

\newcounter{dctr}[section]
\numberwithin{equation}{section}

\newtheorem{deff}{Definition}[section]
\newtheorem{proposition}{Proposition}[section]
\newtheorem{theorem}{Theorem}[section]
\newtheorem{remark}{Remark}[section]

\newcommand{\mig}{\frac{1}{2}}
\newcommand{\bigO}{\mathcal{O}}

\newcommand{\C}{\mathcal{C}}

\renewcommand{\L}{\mathcal{L}}
\def\sii{\Leftrightarrow}

\newcommand{\Z}{\mathbb{Z}}
\newcommand{\N}{\mathbb{N}}

\allowdisplaybreaks
\numberwithin{equation}{section}

\title[Efficient WENO schemes for {nonuniform} grids]{Efficient WENO schemes for {nonuniform} grids}
\author[Marti]{M. Carmen Mart\'{i}$^{\mathrm{a}}$}
\author[Mulet]{Pep Mulet$^{\mathrm{a}}$}
\author[Y\'a\~nez]{Dionisio F. Y\'a\~nez$^{\mathrm{a}}$}
\author[Zorío]{David Zorío$^{\mathrm{a}}$}

\usepackage[normalem]{ulem}

\begin{document}

\begin{abstract}
A set of arbitrarily high-order WENO schemes for reconstructions on {nonuniform} grids is presented. These non-linear interpolation methods use simple smoothness indicators with a linear cost with respect to the order, making them easy to implement and computationally efficient. The theoretical analysis to verify the accuracy and the essentially non-oscillatory properties are presented together with some numerical experiments involving algebraic problems in order to validate them. Also, these general schemes are applied for the solution of conservation laws and hyperbolic systems in the context of finite volume methods.
\end{abstract}

\date{April 25, 2024.}


\thanks{$^{\mathrm{A}}$Departament de Matem\`{a}tiques,
Universitat de Val\`{e}ncia, Av.\ Vicent Andr\'es Estell\'es, E-46100
Burjassot,
 Spain.  E-Mails:
   {\tt Maria.C.Marti@uv.es}, {\tt pep.mulet@uv.es}, {\tt dionisio.yanez@uv.es}, {\tt david.zorio@uv.es}}

\maketitle

\section{Introduction}

\subsection{Scope} \label{subsec:scope}

Non-linear interpolation methods as essentially non-oscillatory (ENO) and weighted ENO (WENO) schemes were developed to approximate numerical solutions of partial differential equations (PDEs), see e.g. \cite{ShuOsher89, ShuOsher1989, JiangShu96, LiuOsherChan94, shu09, zhangshu16}, although they have also been used successfully in non-PDE applications as image processing (see e.g. \cite{ARSY, ABM}), computer vision (see e.g. \cite{sikishu09} mentioned by \cite{shu09}) or finance (see e.g. \cite{KEG}). Recently, a vast quantity of new WENO variants have been proposed to enhance the accuracy of the traditional WENO scheme proposed in \cite{ShuOsher89, ShuOsher1989} as, for instance, Mapped WENO (WENO-M) by Henrick et al. \cite{HenrickAslamPowers2005}, WENO-Z \cite{borgescarmona} or a WENO scheme with progressive order of accuracy close to discontinuities \cite{ARSY20}.

The idea behind the design of the WENO algorithm is quite powerful. It is based on the construction of a piecewise interpolator trying to avoid the points where the function is discontinuous. For this, a convex combination of interpolations of lower degree is proposed, using non-linear weights {for} the contribution of each interpolator. The final result is an interpolator of high-order accuracy if the data used are free of discontinuities and a lower order of accuracy interpolator if the {utilized stencil} crosses a singularity. The main ingredients of WENO algorithms are the optimal weights, which are the values needed to obtain the highest order of accuracy, and the smoothness indicators which indicate if a subset of points contains a discontinuity. There are certain situations where the optimal weights are negative and non-linear WENO weights do not satisfy the required conditions (see \cite{shu09}). Some solutions for this issue are summarized in \cite{shu09}. Levy et al. proposed in \cite{LevyPuppoRusso} a central scheme based on a centered version of WENO, named CWENO, which overcomes the mentioned matter.

On the one hand, in \cite{BaezaBurgerMuletZorio2019}, the authors introduced a central WENO scheme based on a global average weight independent of the optimal weights. This method has two key ingredients: first, the global average weight, which is defined using the smoothness indicators through elementary operations, and second, the evaluation of the interpolator polynomial using the whole stencil. This method was developed for structured grids. However, in some applications (see \cite{CFR05, shihushu02},  mentioned in \cite{shu09}), WENO reconstructions on non-uniform meshes are needed.

On the other hand, in \cite{BBMZ18}, the authors presented an alternative family of smoothness indicators that are simpler and computationally cheaper than those originally proposed by Jiang and Shu in \cite{JiangShu96}. They proved that the computational time needed to obtain the proposed smoothness indicators increased linearly with respect to the accuracy order of the WENO scheme considered while for Jiang and Shu's smoothness indicators the growth is quadratic.

In this paper, we combine and extend the ideas proposed in \cite{BaezaBurgerMuletZorio2019, BBMZ18} for {nonuniform} grids and prove theoretically that the new proposed scheme satisfies the desired accuracy. In addition, we show that this scheme efficiently solves, in some cases, the loss of order of accuracy near smooth extrema associated with the original WENO and CWENO schemes. We will present some mimetic examples to check these theoretical results. Another key feature of this method is that the cost of the computation of the smoothness indicators for a nonuniform stencil does not have a significant impact on the whole computational cost of the scheme, compared with the cost that the computation of the smoothness indicators has on the original WENO and CWENO schemes. We will apply the new WENO reconstruction method in finite volume schemes to obtain numerical solutions in shock problems from hyperbolic conservation laws to illustrate its performance.

\subsection{Outline of this paper}

The remainder of the paper is divided as follows: in Section \ref{secpreliminares} we briefly review the classical Lagrange interpolation and the notation that will be used throughout this paper. We present linear reconstructions from point values and cell averages on {nonuniform} grids. The new WENO method will be introduced in Section \ref{generalizado} and some theoretical results will be provided. In particular, we will prove that the new WENO scheme satisfies the necessary properties to achieve maximum order of accuracy when the data do not cross any singularity and the resulting interpolator is essentially non-oscillatory. Section \ref{summary_algorithms} is devoted to {showing} the computational aspects of the algorithm. In Section \ref{numerexp}, we present some numerical experiments to validate the theoretical results and to test the performance of the new scheme in shock problems from hyperbolic conservation laws. Finally, some conclusions and remarks will be drawn in Section \ref{conclusions}.

\section{Preliminaries}\label{secpreliminares}

We start introducing basic definitions which will be used throughout this work. We also review the definition of linear pointwise reconstructions on {nonuniform} grids from given point values and cell averages.

\begin{deff}
  In the remainder of the paper, we  denote by $\Pi_k$, $k \in \mathbb{N}_0$, the space of polynomials of {\em maximal} degree~$k$, i.e.
  $$\Pi_k=\left\{p(x)=\sum_{i=0}^k a_ix^i: a_i\in\mathbb{R}, \,0\leq i\leq k\right\},$$
  and by $\bar{\Pi}_k$ the space of polynomials of {\em exact} degree~$k$, $k \in \mathbb{N}_0$, i.e.
  $$\bar{\Pi}_k=\left\{p(x)=\sum_{i=0}^k a_ix^i: a_i\in\mathbb{R}, \,0\leq i\leq k-1,\,a_k\neq 0\right\}.$$
\end{deff}

\begin{deff}
  For $\alpha\in\Z$:
\begin{align*}
  f(h)&=\bigO(h^\alpha)\sii \limsup_{h\to 0} \left|\frac{f(h)}{h^\alpha}\right| <\infty\\
  f(h)&=\overline{\bigO}(h^\alpha)\sii \limsup_{h\to 0} \left|\frac{f(h)}{h^\alpha}\right|  <\infty \land \liminf_{h\to 0} \left|\frac{f(h)}{h^\alpha}\right| > 0
\end{align*}
Since, for positive $f, g$,
\begin{align*}
  \limsup_{h\to 0} f(h)g(h)&\leq \limsup_{h\to 0} f(h) \limsup_{h\to
    0}  g(h)\\
 \liminf_{h\to 0} f(h)g(h)&\geq \liminf_{h\to 0} f(h) \liminf_{h\to 0} g(h)
\end{align*}
it follows that
\begin{align*}
\bigO(h^{\alpha}) \bigO(h^{\beta}) &= \bigO(h^{\alpha+\beta}) \\
\overline\bigO(h^{\alpha}) \overline\bigO(h^{\beta}) &=
\overline\bigO(h^{\alpha+\beta}).\\
\end{align*}
\end{deff}

\subsection{Reconstructions on {nonuniform} grids}\label{reconstructionlagrange}

Since we want to design WENO schemes working on {nonuniform} grids, the first step is to describe the procedure to obtain a reconstruction from given data at the desired order. We will focus both on the pointwise reconstructions from point values (the traditional Lagrange interpolation/extrapolation) and the reconstructions from cell averages.

\subsubsection*{Reconstructions from point values}

Let $S=\{x_0,\ldots,x_{R-1}\}$ be a stencil, $x_i<x_j$, $i<j$, and $f_i:=f(x_i)$, $0\leq i\leq R-1$ for some function $f$. It is well-known that the polynomial $p\in\Pi_{R-1}$ satisfying $p(x_i)=f_i$ is given by the Lagrange form:
$$p(x)=\sum_{i=0}^{R-1}\left[\prod_{j=0,j\neq i}^{R-1}\frac{(x-x_j)}{(x_i-x_j)}\right]f_i.$$

\subsubsection*{Reconstructions from cell averages}

Now, let us consider a stencil consisting of cell interfaces given by $S=\{x_0,\ldots,x_R\}$, $x_i<x_j$, $i<j$, being the associated cells $C_i=[x_i,x_{i+1}]$, $0\leq i\leq R-1$ and $\bar{f}_i$ the known values given by
$$\bar{f}_i:=\frac{1}{x_{i+1}-x_i}\int_{x_i}^{x_{i+1}}f(x)dx.$$
Then, let us assume that we want to find $p\in\Pi_{R-1}$ such that $p(x)$ interpolates $f(x)$. To do so, we impose
\begin{equation}\label{pcell}
  \frac{1}{x_{j+1}-x_j}\int_{x_j}^{x_{j+1}}p(x)dx=\bar{f}_j,\quad0\leq j\leq R-1.
\end{equation}

To find $p$, we use the strategy proposed in \cite{SINUM2011} for structured grids, which we next generalize to the {nonuniform} case. Let us consider
\begin{equation}\label{pcell_lagrange}
  p(x)=\sum_{i=0}^{R-1}p_i(x)\bar{f}_i,
\end{equation}
with $p_i\in\Pi_{R-1}$ to be determined, $0\leq i\leq R-1$. Then, imposing \eqref{pcell} on \eqref{pcell_lagrange} and using linearity, we have for a fixed $j$,
\begin{equation}\label{pcell_lagrange_int}
  \frac{1}{x_{i+1}-x_i}\sum_{i=0}^{R-1}\left[\int_{x_j}^{x_{j+1}}p_i(x)dx\right]\bar{f}_i=\bar{f}_j,\quad0\leq j\leq R-1.
\end{equation}
Therefore, \eqref{pcell_lagrange_int} is satisfied if we impose
\begin{equation}\label{pcell_lagrange_cond}
  \int_{x_j}^{x_{j+1}}p_i(x)dx=\delta_{i,j}(x_{i+1}-x_i),\quad0\leq i,j\leq R-1,
\end{equation}
with $\delta_{i,j}$ the Kronecker delta.

Let $P_i\in\Pi_R$ be a primitive of $p_i$, namely, such that $P_i'=p_i$. Then, for each $0\leq i\leq R-1$, \eqref{pcell_lagrange_cond} reads
\begin{equation*}
  P_i(x_{j+1})-P_i(x_j)=\delta_{i,j}(x_{i+1}-x_i),\quad0\leq j\leq R-1.
\end{equation*}

Hence, it suffices to impose
\[
P_i(x_j)=\begin{cases}
0 & 0\leq j\leq i, \\
x_{i+1}-x_i & i+1\leq j\leq R,
\end{cases}
\]
to write $P_i$ in Lagrange form as
\begin{equation}\label{prim_lagrange}
  P_i(x)=(x_{i+1}-x_i)\sum_{j=i+1}^R\prod_{k=0,k\neq j}^R\frac{x-x_k}{x_j-x_k}.
\end{equation}

Thus, taking derivatives at both sides of \eqref{prim_lagrange} we get
\begin{equation}\label{lagrange}
  p_i(x)=(x_{i+1}-x_i)\sum_{j=i+1}^R\sum_{k=0,k\neq j}^R\frac{1}{x_j-x_k}\prod_{l=0,l\neq j,k}^R\frac{x-x_l}{x_k-x_l}
\end{equation}
and \eqref{pcell_lagrange} then reads
\begin{equation}\label{pcell_final}
p(x)=\sum_{i=0}^{R-1}(x_{i+1}-x_i)\left[\sum_{j=i+1}^R\sum_{k=0,k\neq j}^R\frac{1}{x_j-x_k}\prod_{l=0,l\neq j,k}^R\frac{x-x_l}{x_k-x_l}\right]\bar{f}_i.
\end{equation}

\section{Generalized WENO schemes}\label{generalizado}

The classical Lagrange methods introduced in Section \ref{reconstructionlagrange} are useful to interpolate smooth data but some artifacts appear when they are used in applications with non-continuous inputs. To solve this problem, we present in this section a generalization of the WENO method for {nonuniform} grids.

\subsection{Construction of the method and theoretical results} 

Let $h>0$, $R\geq3$ and let $S_{R,h}=\{x_{0,h},\ldots,x_{R-1,h}\}$ be a stencil, with $x_{i,h}=z+c_ih$, for some fixed $z\in\R$ and $c_i\in\R$, $0\leq i\leq R-1$, $c_i<c_j$, $i<j$. Let $f$ be a function and define $f_{i,h}=f(x_{i,h})$.

We denote by $x_h^*\in\R$ the point in which we desire to reconstruct the data from the stencil $S_{R,h}$.
%
This point is located between $x_{r-1,h}$ and $x_{r,h}$ if $R$ is even, with $r=R/2$, or rather between $x_{r-1,h}$ and $x_{r,h}$ (left-biased) or between $x_{r,h}$ and $x_{r+1,h}$ (right-biased) if $R$ is odd, with $r=(R-1)/2$. Note that this is actually a generalization of the traditional centered or left/right-biased WENO/CWENO schemes on cartesian grids to {nonuniform} grids, respectively, (see e.g. \cite{LiuOsherChan94, LevyPuppoRusso}).

Let us first consider $p_{R,i,h}$ the reconstruction polynomials satisfying $p_{R,i,h}(x_{j,h})=\L[f](x_{j,h})$, where
the operator $\L$ is defined depending on the framework chosen (point values or cell average reconstructions), i.e.,
$$\L[f](x_{j,h})=
\begin{cases}
f(x_{j,h}), &\text{for point value reconstructions}, \\
\bar{f}_{j,h},& \text{for cell average reconstructions}, \\
\end{cases}
$$
with $i\leq j\leq i+r$, for $0\leq i\leq r'$, with $r=\lfloor (R-1)/2\rfloor$ and $r'=\lceil (R-1)/2\rceil$, and evaluate them at $x_h^*$, where $\lfloor \cdot \rfloor$ denotes the floor function and $\lceil \cdot \rceil$ the ceil function.

Following the same idea based on the Yamaleev-Carpenter weight construction proposed in \cite{YamaleevCarpenter2009}, we define $d_{R,h}$ as
$$d_{R,h}=(f[x_{0,h},\ldots,x_{R-1,h}])^2,$$
where $f[x_{0,h},\ldots,x_{R-1,h}]$ is the undivided difference defined as
\begin{align*}
  f[x_{i,h}]&=f_{i,h},\quad 0\leq i\leq R, \\
  f[x_{i,h},\ldots,x_{i+j,h}]&=\frac{f[x_{i+1,h},\ldots,x_{i+j,h}]-f[x_{i,h},\ldots,x_{i+j-1,h}]}{c_{i+j}-c_i},\begin{array}{c}\quad 0\leq i\leq R-j,\\\quad 0<j\leq R.\end{array}
\end{align*}

The simplified smoothness indicators, akin to the idea proposed in \cite{BBMZ18}, which in turn is the key to the simplicity of this algorithm, are defined by
\begin{equation}\label{eq:SI}
I_{R,i,h}:=\sum_{j=i}^{r+i-1}\left(\frac{f_{j+1,h}-f_{j,h}}{c_{j+1}-c_j}\right)^2,\quad0\leq i\leq r'.
\end{equation}

Since in this case, in the context of an arbitrary grid, there is no point in computing ideal weights in order to build the reconstruction with optimal accuracy from the lower order reconstruction, we simply consider a uniform convex combination based on the $r'+1$ values $\frac{1}{r'+1}$ and define
$$\omega_{R,i,h}=\frac{\alpha_{R,i,h}}{\sum\limits_{j=0}^{r'}\alpha_{R,j,h}},\quad\textnormal{with }\,\alpha_{R,i,h}=\frac{1}{r'+1}\left(1+\frac{d_{R,h}^s}{I_{R,i,h}^s+\varepsilon}\right),\quad 0\leq i\leq r',$$
with $s$ an exponent to be determined (will be done so when the accuracy analysis is performed) and $\varepsilon>0$ a small positive quantity in order to avoid divisions by zero. With this weight design, we can define the following weighted essentially non-oscillatory reconstruction:
$$\tilde q_{R,h}(x_h^*)=\sum_{i=0}^{r'}\omega_{R,i,h}p_{R,i,h}(x_h^*).$$

Since this reconstruction does not attain the desired optimal order, $R$, we use the reconstruction associated with the whole stencil $S_{R,h}$, $p_{R,h}$ satisfying $p_{R,h}(x_{i,h})=\L[f](x_{i,h})$, which we evaluate at $x_h^*$, and the following global average weight, given by
$$\omega_{R,h}=\frac{1}{1+d_{R,h}^sJ_{R,h,s}},\quad\textnormal{with }\,J_{R,h,s}=\sum_{i=0}^{r'}\frac{1}{\displaystyle I_{R,i,h}^s+\varepsilon}.$$

Then, the final reconstruction with essentially non-oscillatory properties is given by
$$q_{R,h}(x_h^*)=\omega_{R,h}p_{R,h}(x_h^*)+(1-\omega_{R,h})\tilde q_{R,h}(x_h^*).$$

Now that we have introduced and defined all the components of the WENO method, we are in a position to prove the principal properties of the scheme in terms of the order of the weights and accuracy of the reconstructions obtained.

\begin{proposition}\label{weights_standard}
  For the standard case, and keeping the same notation, if $f$ has a critical point of order $k$ at $z$, $0\leq k\leq R-2$, then it holds $0\leq\omega_{R,i,h}\leq1$, $0\leq\omega_{R,h}\leq1$ and
  \begin{align*}
    \omega_{R,i,h}&=\begin{cases}
    \displaystyle\frac{1}{r'+1}+\bigO(h^{2s(R-1-m)})+\bigO(\varepsilon) & \textnormal{if }S_{R,h}\textnormal{ is smooth,}\\
    \bigO(1)+\bigO(\varepsilon) & \textnormal{if a discontinuity crosses }S_{R,h}\setminus S_{R,i,h},\\
    \bigO(h^{2s})+\bigO(\varepsilon) & \textnormal{if a discontinuity crosses }S_{R,i,h},
    \end{cases}\\
    \omega_{R,h}&=\begin{cases}
    1{+}\bigO(h^{2s(R-1-m)}){+}\bigO(\varepsilon) & \textnormal{if }S_{R,h}\textnormal{ is smooth,}\\
    \bigO(h^{2s})+\bigO(\varepsilon) & \textnormal{if a discontinuity crosses }S_{R,h},
    \end{cases}
  \end{align*}
  with $m=\max\limits_{0\leq i<R-1}m_i$, where
  \[m_i=
  \begin{cases}
    \min\{q\in\N\colon 2 | q, q\geq k, f^{(q+1)}(z)\neq0\}+1 & \textnormal{if }3\leq R\leq 4\textnormal{ and }c_{i}+c_{i+1}=0,\\
    k+1 & \textnormal{otherwise}.
  \end{cases}
  \]
\end{proposition}

\begin{proof}
  If $f\in\C^R$, then $d_{R,h}=\bigO(h^{2R-2})$. On the other hand, by all the considerations derived from \cite[Section 2]{BaezaBurgerMuletZorio2018} and \cite[Section 2]{BBMZ18} we conclude that, since $r=1$ if $3\leq R\leq 4$ and $r>1$ if $R\geq5$, it holds $I_{R,i,h}=\overline\bigO(h^{2m_i})$, with
  \[
  m_i=
  \begin{cases}
    \min\{q\in\N\colon 2 | q, q\geq k, f^{(q+1)}(z)\neq0\}+1 & \textnormal{if }3\leq R\leq 4\textnormal{ and }c_i+c_{i+1}=0,\\
    k+1 & \textnormal{otherwise}.
  \end{cases}
  \]
Therefore, we have
\begin{align*}
  \alpha_{R,i,h}&=\frac{1}{r'+1}\left(1+\frac{d_{R,h}^s}{I_{R,i,h}^s+\varepsilon}\right)=\frac{1}{r'+1}\left(1+\frac{\bigO(h^{2s(R-1)})}{\overline\bigO(h^{2sm_i})}\right)+\bigO(\varepsilon)\\&=\frac{1}{r'+1}+\bigO(h^{2s(R-1-m_i)})+\bigO(\varepsilon).
\end{align*}
Hence,
\begin{align*}
  \omega_{R,i,h}&=\frac{\alpha_{R,i,h}}{\sum\limits_{j=0}^{r'}\alpha_{R,j,h}}=\frac{\frac{1}{r'+1}+\bigO(h^{2s(R-1-m_i)})+\bigO(\varepsilon)}{\sum\limits_{j=0}^{r'}\left(\frac{1}{r'+1}+\bigO(h^{2s(R-1-m_j)})+\bigO(\varepsilon)\right)}\\
  &=\frac{\frac{1}{r'+1}+\bigO(h^{2s(R-1-m_i)})+\bigO(\varepsilon)}{1+\bigO(h^{2s(R-1-m)})+\bigO(\varepsilon)}=\frac{1}{r'+1}+\bigO(h^{2s(R-1-m)})+\bigO(\varepsilon),
\end{align*}
with $m=\max\limits_{0\leq i<R-1}m_i$.

Now, let us assume that a discontinuity crosses $S_{R,h}$. Then, we define as $A$ the set of indexes $j$ such that the discontinuity does not cross $S_{R,{j},h}$ and as $B$ the set of indexes $j$ such that the discontinuity crosses $S_{R,{j},h}$. Note that, by construction of $r'$, the stencils {do} not overlap, and thus $A\neq\emptyset$. Then, on the one hand, if $j\in A$:
\begin{align*}
  \alpha_{R,j,h}&=\frac{1}{r'+1}\left(1+\frac{d_{R,h}^s}{I_{R,j,h}^s+\varepsilon}\right)=\frac{1}{r'+1}\left(\frac{\overline\bigO(1)}{\overline\bigO(h^{2sm_j})}\right)+\bigO(\varepsilon)=\overline\bigO(h^{-2sm_j})+\bigO(\varepsilon).
\end{align*}
On the other hand, if $j\in B$:
\begin{align*}
  \alpha_{R,j,h}&=\frac{1}{r'+1}\left(1+\frac{d_{R,h}^s}{I_{R,j,h}^s+\varepsilon}\right)=\frac{1}{r'+1}\left(\frac{\overline\bigO(1)}{\overline\bigO(1)}\right)+\bigO(\varepsilon)=\overline\bigO(1)+\bigO(\varepsilon).
\end{align*}

Consequently, if $i\in A$, we have:
\begin{align*}
  \omega_{R,i,h}&=\frac{\alpha_{R,i,h}}{\sum\limits_{j=0}^{r'}\alpha_{R,j,h}}=\frac{\alpha_{R,i,h}}{\sum\limits_{j\in A}\alpha_{R,j,h}+\sum\limits_{j\in B}\alpha_{R,j,h}}\\
  &=\frac{\overline\bigO(h^{-2sm_i})+\bigO(\varepsilon)}{\sum\limits_{j\in A}(\overline\bigO(h^{-2sm_j})+\bigO(\varepsilon))+\sum\limits_{j\in B}(\overline\bigO(1)+\bigO(\varepsilon))}=\frac{\overline\bigO(h^{-2sm_i})+\bigO(\varepsilon)}{\overline\bigO(h^{-2sm'})+\bigO(\varepsilon)}\\&=\bigO(1)+{\bigO(\varepsilon)},
\end{align*}
with $m':=\max\limits_{j\in A}m_j$, which satisfies $m'\geq m_i$, since $i\in A$ by assumption. However, if $i\in B$, then:
\begin{align*}
  \omega_{R,i,h}&=\frac{\alpha_{R,i,h}}{\sum\limits_{j=0}^{r'}\alpha_{R,j,h}}=\frac{\alpha_{R,i,h}}{\sum\limits_{j\in A}\alpha_{R,j,h}+\sum\limits_{j\in B}\alpha_{R,j,h}}\\
  &=\frac{\overline\bigO(1)+\bigO(\varepsilon)}{\sum\limits_{j\in A}(\overline\bigO(h^{-2sm_j})+\bigO(\varepsilon))+\sum\limits_{j\in B}(\overline\bigO(1)+\bigO(\varepsilon))}=\frac{\overline\bigO(1)+\bigO(\varepsilon)}{\overline\bigO(h^{-2sm'})+\bigO(\varepsilon)}\\&=\bigO(h^{2s})+{\bigO(\varepsilon)}.
\end{align*}

Now we have to analyze $\omega_{R,h}$. To do so, we first study the accuracy of $J_{R,h,s}$:
\begin{align*}
  J_{R,h,s}&=\sum_{i=0}^{r'}\frac{1}{I_{R,i,h}^s+\varepsilon}=\sum_{i\in A}\frac{1}{I_{R,i,h}^s+\varepsilon}+\sum_{i\in B}\frac{1}{I_{R,i,h}^s+\varepsilon}\\
  &=\sum_{i\in A}\frac{1}{\overline\bigO(h^{2sm_i})}+\sum_{i\in B}\frac{1}{\overline\bigO(1)}+\bigO(\varepsilon)=\sum_{i\in A}\overline\bigO(h^{-2sm_i})+\sum_{i\in B}\overline\bigO(1)+\bigO(\varepsilon)\\
  &=\overline\bigO(h^{-2sm'})+\bigO(\varepsilon).
\end{align*}
Therefore, if $f\in\C^R$, then $d_{R,h}=\bigO(h^{2R-2})$ and
\begin{align*}
  \omega_{R,h}&=\frac{1}{1+d_{R,h}J_{R,h,s}}=\frac{1}{1+\bigO(h^{2R-2})(\overline\bigO(h^{-2sm'})+\bigO(\varepsilon))}\\
  &=1{+}\bigO(h^{-2s(R-1-m')}){+}\bigO(\varepsilon)=1{+}\bigO(h^{-2s(R-1-m)}){+}\bigO(\varepsilon),
\end{align*}
where the last equality holds since $m'\leq m$. Finally, if a discontinuity crosses $S_{R,h}$, then $d_{R,h}=\overline\bigO(1)$ and
\begin{align*}
  \omega_{R,h}&=\frac{1}{1+d_{R,h}J_{R,h,s}}=\frac{1}{1+\overline\bigO(1)(\overline\bigO(h^{-2sm'})+\bigO(\varepsilon))}=\bigO(h^{2sm'})+\bigO(\varepsilon)\\&=\bigO(h^{2s})+\bigO(\varepsilon).
\end{align*}
\end{proof}

\begin{theorem}\label{acc_standard}
  Keeping the same notation, if $s\geq (r+1)/2$, $z$ is not a critical point of $f$ for $3\leq R\leq4$ or $z$ is not a critical point of $f$ of order $R-2$ for $R\geq5$, then the optimal accuracy is attained both in the case of smoothness and the case of discontinuities, namely:
  \[
  q_{R,h}(x_h^*)=
  \begin{cases}
    f(x_h^*)+\bigO(h^R)+\bigO(\varepsilon) & f\in\C^R, \\
    f(x_h^*)+\bigO(h^{r+1})+\bigO(\varepsilon) & \textnormal{if a discontinuity crosses }S_{R,h}.
  \end{cases}
  \]
\end{theorem}

\begin{proof}
  Let us first assume that $f\in\C^R$. Then, by Proposition \ref{weights_standard}:
  \begin{equation}\label{eq:bounds}
    \begin{split}
      \omega_{R,i,h}&=\displaystyle\frac{1}{r'+1}+\bigO(h^{2s(R-1-m)})+\bigO(\varepsilon),\\
      \omega_{R,h}&=1{+}\bigO(h^{2s(R-1-m)}){+}\bigO(\varepsilon),
    \end{split}
  \end{equation}
  with $m=\max\limits_{0\leq i<R-1}m_i$, where
  \[m_i=
  \begin{cases}
    \min\{q\in\N\colon 2 | q, q\geq k, f^{(q+1)}(z)\neq0\}+1 & \textnormal{if }3\leq R\leq 4\textnormal{ and }c_{i}+c_{i+1}=0,\\
    k+1 & \textnormal{otherwise}.
  \end{cases}
  \]

  Now, taking into account our hypothesis involving the critical points, we have that $m_i=k+1$, $\forall\, 0\leq i\leq R-1$, $\forall\, R\geq3$. Moreover, $k=0$ for $3\leq R\leq 4$ and $k\neq R-2$ for $R\geq5$, and thus we can assume $k<R-2$ for any $R\geq3$ (the case $k>R-2$ being trivial, as we will point out at the end of the proof). Therefore, taking into account that
  \begin{equation}\label{eq:acc}
    \begin{split}
      p_{R,i,h}(x_h^*)&=f(x_h^*)+\bigO(h^{r+1}),\\
      p_{R,h}(x_h^*)&=f(x_h^*)+\bigO(h^R),
    \end{split}
  \end{equation}
  and using that $\sum\limits_{i=0}^{r'}\omega_i=1$, we have
  \begin{align*}
    \tilde q_{R,h}(x_h^*)&=\sum_{i=0}^{r'}\omega_ip_{R,i,h}(x_h^*)=\sum_{i=0}^{r'}\left[\omega_i\left(f(x_h^*)+\bigO(h^{r+1})\right)\right]=\sum_{i=0}^{r'}\omega_if(x_h^*)+\bigO(h^{r+1})\\
    &=f(x_h^*)+\bigO(h^{r+1}).
  \end{align*}

Using \eqref{eq:bounds} and \eqref{eq:acc}:
  \begin{align*}
    q_{R,h}(x_h^*)&=\omega_{R,h}p_{R,h}(x_h^*)+(1-\omega_{R,h})\tilde q_{R,h}(x_h^*)\\
    &=\omega_{R,h}\left(f(x_h^*)+\bigO(h^R)\right)+(1-\omega_{R,h})\left(f(x_h^*)+\bigO(h^{r+1})\right)\\
    &=f(x_h^*)+\omega_{R,h}\bigO(h^R)+(1-\omega_{R,h})\bigO(h^{r+1})\\
    &=f(x_h^*)+\left(1{+}\bigO(h^{2s(R-2-k)}){+}\bigO(\varepsilon)\right)\bigO(h^R)\\
    &+\left(\bigO(h^{2s(R-2-k)})+\bigO(\varepsilon)\right)\bigO(h^{r+1})=f(x_h^*)+\bigO(h^{2s(R-2-k)+r+1})+\bigO(\varepsilon),
  \end{align*}
 where in the last equality we have used that $R\geq r+1$. Therefore, if one wants the optimal $R$-th order accuracy, one must impose
  \begin{align*}
    2s(R-2-k)+r+1\geq R\sii s\geq\frac{R-r-1}{2(R-2-k)}\geq\frac{R-r-1}{2}\geq\frac{r}{2}.
  \end{align*}
  On the other hand, since
  $$\frac{R-r-1}{2(R-2-k)}\leq\frac{R-r-1}{2}\leq\frac{r+1}{2},$$
  it suffices to impose
  $$s\geq\frac{r+1}{2}.$$

Let us consider now the case where a discontinuity crosses $S_{R,h}$, then, again by Proposition \ref{weights_standard}, we have
  \begin{equation}\label{eq:bounds_disc}
    \begin{split}
      \omega_{R,i,h}&=\begin{cases}
      \displaystyle\bigO(1)+\bigO(\varepsilon) & i\in A,\\
      \displaystyle\bigO(h^{2s})+\bigO(\varepsilon) & i\in B,
      \end{cases}\\
      \omega_{R,h}&=\bigO(h^{2s})+\bigO(\varepsilon),
    \end{split}
  \end{equation}
  and now it holds
  \begin{equation}\label{eq:acc_disc}
    \begin{split}
      p_{R,i,h}(x_h^*)&=\begin{cases}
      f(x_h^*)+\bigO(h^{r+1}) & i\in A,\\
      \bigO(1) & i\in B,
      \end{cases}\\
      p_{R,h}(x_h^*)&=\bigO(1).
    \end{split}
  \end{equation}
  Therefore,
  \begin{align*}
    \tilde q_{R,h}(x_h^*)&=\sum_{i=0}^{r'}\omega_ip_{R,i,h}(x_h^*)=\sum_{i\in A}\omega_ip_{R,i,h}(x_h^*)+\sum_{i\in B}\omega_ip_{R,i,h}(x_h^*)\\
    &=\sum_{i\in A}\omega_i(f(x_h^*)+\bigO(h^{r+1}))+\sum_{i\in B}{(}\bigO(h^{2s}){+\bigO(\varepsilon))}\bigO(1)\\
    &=f(x_h^*)\sum_{i\in A}\omega_i+\bigO(h^{2s})+\bigO(h^{r+1}){+\bigO(\varepsilon)}\\
    &=f(x_h^*)\left(1-\sum_{i\in B}\omega_i\right)+\bigO(h^{2s})+\bigO(h^{r+1}){+\bigO(\varepsilon)}\\
    &=f(x_h^*)\left(1-\sum_{i\in B}\bigO(h^{2s})\right)+\bigO(h^{2s})+\bigO(h^{r+1}){+\bigO(\varepsilon)}\\
    &=f(x_h^*)\left(1{+}\bigO(h^{2s})\right)+\bigO(h^{2s})+\bigO(h^{r+1})=f(x_h^*)+\bigO(h^{2s})+\bigO(h^{r+1}){+\bigO(\varepsilon)}.
  \end{align*}
Combining in this case \eqref{eq:bounds_disc} and \eqref{eq:acc_disc}, we can conclude that
  \begin{align*}
    q_{R,h}(x_h^*)&=\omega_{R,h}p_{R,h}(x_h^*)+(1-\omega_{R,h})\tilde q_{R,h}(x_h^*)\\
    &=\left(\bigO(h^{2s})+\bigO(\varepsilon)\right)\bigO(1)+\left(1{+}\bigO(h^{2s}){+}\bigO(\varepsilon)\right)\left(f(x_h^*)+\bigO(h^{2s})+\bigO(h^{r+1}){+\bigO(\varepsilon)}\right)\\
    &=f(x_h^*)+\bigO(h^{2s})+\bigO(h^{r+1})+\bigO(\varepsilon).
  \end{align*}
  Then, in order to attain the optimal accuracy, one must impose
  $$2s\geq r+1\sii s\geq\frac{r+1}{2},$$
  and therefore, it suffices with the choice
  $$s\geq\frac{r+1}{2}$$
  to attain optimal accuracy both in the smooth and the discontinuous case.
\end{proof}

\subsection{Remarks regarding the computational cost}\label{sec:comp}
 
As stated in the introduction, one of the strongest points of the proposed WENO method is the ease of computation of the smoothness indicators, even taking into consideration the fact that it is defined in the context of a nonuniform grid.

Unlike the case of the traditional Jiang-Shu smoothness indicators, which in turn are used in the method proposed in \cite{BaezaBurgerMuletZorio2019}, the indicators used in this work, based on those proposed in \cite{BBMZ18}, are considerably cheaper to compute, especially when the grid is non-uniform.

To be more specific, given $R\in\N$, a stencil of the form $S_{R,h}=\{x_{0,h},\ldots,x_{R-1,h}\}$, with $x_{i,h}=z+c_ih$, $c_i<c_j$, the set of nodal values $\{f_{0,h},\ldots,f_{R-1,h}\}$, with $f_{i,h}=f(x_{i,h})$, for $f$ a function to be approximated, and the reconstruction point $x_h^*=z+c^*h$, the simplified expression for the classical Jiang-Shu smoothness indicators has the structure
\begin{equation}\label{eq:JS}
\sum_{i=0}^{R-1}\sum_{j=i}^{R-1}\left(\sum_{k=0}^{R-1}\alpha_{i,j,k}c_k\right)f_{i,h}f_{j,h},
\end{equation}
with $\alpha_{i,j,k}$ to be determined for each $0\leq i\leq j\leq R-1$, $0\leq k\leq R-1$.

From \eqref{eq:JS}, it can be clearly seen that the number of operations that are required to compute the Jiang-Shu smoothness indicators in terms of the accuracy order is $\bigO(R^3)$. Therefore, any WENO method using these smoothness indicators has, at least, a cubic cost with respect to the aforementioned order of the method. This includes the techniques proposed by Levy, Puppo and Russo in 
\cite{LevyPuppoRusso} or those proposed in \cite{BaezaBurgerMuletZorio2019}.

In contrast, from \eqref{eq:SI} the smoothness indicators proposed in this paper have linear cost with respect to the accuracy order, namely, $\bigO(R)$. It is also important to mention that, depending on the coefficients $c_i$ and the location of the reconstruction point $x_h^*$, the ideal weights in the method proposed by Levy, Puppo and Russo can be negative. This issue was fixed in \cite{BaezaBurgerMuletZorio2019} and the method proposed here neither suffers from this phenomenon.

Since computing an interpolating polynomial under a non-uniform grid requires evaluating an expression of the form
\begin{equation}\label{eq:IP}
\sum_{i=0}^{R-1}\left(\sum_{k=0}^{R-1}\beta_{i,k}c_k\right)f_{i,h},
\end{equation}
with $\beta_{i,k}$ to be determined, $0\leq i,k\leq R-1$, the cost from that part of the algorithm requires a quadratic number of operations with respect to the accuracy order, namely, $\bigO(R^2)$. Thus, taking into account that $\lfloor\frac{R}{2}\rfloor$ polynomials and smoothness indicators have to be computed in the process, the computational cost of the WENO methods using the Jiang-Shu smoothness indicators have $\bigO(R^4)$ cost when they are applied on nonuniform grids, while the WENO method proposed in this paper has $\bigO(R^3)$ cost with respect to the accuracy order in a context of nonuniform grids.

In the particular case in which a static nonuniform grid is used during the whole time evolution, the terms $\alpha_{i,j,k}$ and $\beta_{i,k}$ can be precomputed for each local stencil and stored to be used on each time iteration. This implies that the computational cost with respect to the order of the expression in \eqref{eq:JS} is reduced to quadratic, while the cost in \eqref{eq:IP} is reduced to linear. Therefore, the global computational cost of the CWENO methods \cite{BaezaBurgerMuletZorio2019, LevyPuppoRusso} can be reduced to $\bigO(R^3)$, while the cost associated with the method presented in this paper is reduced to $\bigO(R^2)$.

\section{Summary of the algorithms}\label{summary_algorithms}

In this Section, we summarize the algorithm of the proposed WENO method, including the considerations involving the items to attain unconditionally the optimal accuracy. We divide this section into two parts: First, we introduce the description for the reconstruction from point values and after we present the necessary changes to adapt it for reconstructions from cell average values.

\subsection{Algorithm for reconstructions from point values}  $\,\,$

Input data:
  \begin{itemize}
    \item Grid of nodes: $S_{R,h}=\{x_{0,h},\ldots,x_{R-1,h}\}$, with $x_{i,h}=z+c_ih$, $c_i<c_j$, $i<j$, $0\leq i,j\leq R-1$. It suffices to provide the $c_i$ values, $0\leq i\leq R-1$.

    \item Nodal values: $\{f_{0,h},\ldots,f_{R-1,h}\}$, with $f_{i,h}=f(x_{i,h})$.

    \item Reconstruction point: $x_h^*$, satisfying $x_{R/2-1,h}\leq x_h^*\leq x_{R/2,h}$ if $R$ is even or $x_{(R-1)/2-1,h}\leq x_h^*\leq x_{(R-1)/2+1,h}$ if $R$ is odd. It suffices to provide $c^*$, with $x_h^*=z+c^*h$.
    
    \item Small positive quantity: $\varepsilon>0$.
  \end{itemize}

Then, the next steps are to be followed:
\begin{enumerate}
\item Compute
  \begin{align*}
  r&=\left\lfloor\frac{R-1}{2}\right\rfloor,\\
  r'&=\left\lceil\frac{R-1}{2}\right\rceil,\\
  s&=\left\lceil\frac{r+1}{2}\right\rceil.
  \end{align*}

\item Compute the reconstruction polynomials at $c_h^*$ of the corresponding substencils $S_{R,i,h}$, $0\leq i\leq r'$, $p_{R,i,h}$, and the one associated to the global stencil, $S_{R,h}$, $p_{R,h}$, which are:
    \begin{align*}
      p_{R,i,h}(c_h^*)&=\sum_{j=i}^{i+r}\left[\prod_{l=i,l\neq j}^{i+r}\frac{(c_h^*-c_l)}{(c_j-c_l)}\right]f_j,\\
      p_{R,h}(c_h^*)&=\sum_{j=0}^{R-1}\left[\prod_{l=0,l\neq j}^{R-1}\frac{(c_h^*-c_l)}{(c_j-c_l)}\right]f_j.
    \end{align*}

\item Compute the smoothness indicators:
    $$I_{R,i,h}:=\sum_{j=i}^{r+i-1}\left(\frac{f_{j+1,h}-f_{j,h}}{c_{j+1}-c_j}\right)^2,\quad0\leq i\leq r'.$$
  Optionally, the denominators can be omitted in both expressions without any impact on the accuracy of the scheme.

\item Obtain $d_{R,h}$:
      $$d_{R,h}=\left((R-1)!\sum_{i=0}^{R-1}\frac{1}{\displaystyle\prod_{j=0,j\neq i}^{R-1}(c_i-c_j)}f_i\right)^2.$$

\item Compute $J_{R,h,s}$:
  $$J_{R,h,s}=\sum_{i=0}^{r'}\frac{1}{I_{R,i,h}^s+\varepsilon}.$$

\item Generate the weights:
  \begin{align*}
    \omega_{R,i,h}&=\frac{\alpha_{R,i,h}}{\sum\limits_{j=0}^{r'}\alpha_{R,j,h}},\textnormal{ with }\alpha_{R,i,h}=\frac{1}{r'+1}\left(1+\frac{d_{R,h}^s}{I_{R,i,h}^s+\varepsilon}\right),\\
    \omega_{R,h}&=\frac{1}{1+d_{R,h}^sJ_{R,h,s}}.
  \end{align*}

\item Obtain the half-optimal-order essentially non-oscillatory reconstruction:
  $$\tilde q_{R,h}(x_h^*)=\sum_{i=0}^{r'}\omega_ip_{R,i,h}(x_h^*).$$

\item Compute the optimal order essentially non-oscillatory reconstruction:
  $$q_{R,h}(x_h^*)=\omega_{R,h}p_{R,h}(x_h^*)+(1-\omega_{R,h})\tilde q_{R,h}(x_h^*).$$
\end{enumerate}

\subsection{Algorithm for reconstructions from cell average values} $\,$

Input data:
  \begin{itemize}
    \item Grid of cells: $S_{R,h}=\{C_{0,h},\ldots,C_{R-1,h}\}$, with $C_i=[x_{i,h},x_{i+1,h}]$, $0\leq i\leq R-1$, $x_{i,h}=z+c_ih$, $c_i<c_j$, $i<j$, $0\leq i,j\leq R$. It suffices to provide the $c_i$ values, $0\leq i\leq R$.

    \item Cell averages: $\{\bar{f}_{0,h},\ldots,\bar{f}_{R-1,h}\}$, with
      $$\bar{f}_{i,h}=\frac{1}{x_{i+1,h}-x_{i,h}}\int_{x_{i,h}}^{x_{i+1,h}}f(x)dx.$$

    \item Reconstruction point: $x_h^*$, satisfying $x_{(R+1)/2-1,h}\leq x_h^*\leq x_{(R+1)/2,h}$ if $R$ is odd or $x_{R/2-1,h}\leq x_h^*\leq x_{R/2+1,h}$ if $R$ is even. It suffices to provide $c^*$, with $x_h^*=z+c^*h$.
    
    \item Small positive quantity: $\varepsilon>0$.
  \end{itemize}

In step (2) reconstruction polynomials should be defined by:
    \begin{align*}
      p_{R,i,h}(c_h^*)&=\sum_{j=i}^{i+r}(x_{j+1}-x_j)\left[\sum_{l=j+1}^{i+r+1}\sum_{m=i,m\neq l}^{i+r+1}\prod_{n=i,n\neq l,m}^{i+r+1}\frac{c_h^*-c_n}{c_m-c_n}\right]\bar{f}_j,\\
      p_{R,h}(c_h^*)&=\sum_{j=0}^{R-1}(x_{j+1}-x_j)\left[\sum_{l=j+1}^R\sum_{m=0,m\neq l}^R\prod_{n=0,n\neq l,m}^R\frac{c_h^*-c_n}{c_m-c_n}\right]\bar{f}_j.
    \end{align*}

In step (3) the smoothness indicators are:
    $$I_{R,i,h}:=\sum_{j=i}^{r+i-1}\left(\frac{\bar{f}_{i+1,h}-\bar{f}_{i,h}}{\bar{c}_{i+1}-\bar{c}_i}\right)^2,\quad0\leq i\leq r'n$$
and in step (4):
    $$d_{R,h}=\left((R-1)!\sum_{i=0}^{R-1}(x_{i+1}-x_i)\left[\sum_{j=i+1}^R\sum_{k=0,k\neq j}^R\frac{1}{\displaystyle\prod_{l=0,l\neq j,k}^R(c_k-c_l)}\right]\bar{f}_i\right)^2.$$

\begin{remark}
 Note that, if one works with fixed grids, then all the terms involving computations related to the grid spacings, namely, $c_i$, as well as the relative position of the reconstructions points $c_h^*$, can be computed and stored only once, and then used to perform the computations, thus reducing to linear cost with respect to the order the computation of terms such as $p_{R,i,h}$, $p_{R,h}$ and $d_{R,h}$, yielding a comparable performance with respect to the classical WENO schemes with structured grids.
\end{remark}

\section{Numerical experiments}\label{numerexp}

In this section, we perform some numerical examples to check the accuracy properties of our schemes. We divide the experiments into two sections: Firstly, we apply our algorithm in mimetic algebraic numerical tests to assess the order of accuracy. Secondly, we numerically solve several PDEs. In both cases, we implement the method to non-continuous data.

In all the experiments the C++ wrapper \cite{Holoborodko} of the GNU MPFR library \cite{MPFR} is used with a precision of $332$ bits ($\approx 100$ digits) and taking $\varepsilon=10^{-10^5}$.

\subsection{Algebraic numerical experiments}

The goal of this subsection is to check the behaviour of the proposed generalized WENO scheme in different scenarios. We present two different tests: one with a smooth function, and one for non-smooth functions computed both for point values and for cell average values for which the estimated order of accuracy is calculated.

In all the numerical experiments, we consider a mesh of the form $x_{i,h}=z+c_ih_n$, $0\leq i<N$, with $h_{n+1}=\frac{h_n}{2}$ and provide the table of the resulting experiments based on taking $h_0=0.2$ and $0\leq n\leq 20$. The order values are thus computed as
\begin{equation}\label{order}
o_n=\log_2\left(\frac{E_{n-1}}{E_n}\right),\quad E_n=|u_{h_n}^*-f(x_{h_n}^*)|,\quad 0\leq n\leq 20,
\end{equation}
with $x_{h_n}^*$ the reconstruction node and $u_{h_n}^*$ the corresponding resulting WENO reconstruction.

\begin{table}
\begin{center}
        \begin{tabular}{|l|rr|rr|}       \hline
         & \multicolumn{2}{c|}{Test 1} &  \multicolumn{2}{c|}{Test 2} \\
        \hline
        $n$ & Point values & Cell averages & Point values & Cell averages \\
        \hline
         $c_0$  & -3.5411 & -3.5451  & -1.5411 & -3.5451   \\
         $c_1$  & -2.8706 & -2.9810  & -0.9907 & -2.9810\\
         $c_2$  & -2.1411 & -2.3102  & 0.0000 & -2.3102\\
         $c_3$  & -1.7503 & -2.1178  & 0.6792 & -2.1178\\
         $c_4$  & -0.9907 & -1.4574  & 1.7413 & -0.1231\\
         $c_5$  & -0.2145 & -0.8571  & 2.5614 & 0.0000\\
         $c_6$  & 0.6792 & 0.1245    & 3.1410 & 0.8073\\
         $c_7$  & 1.3204 & 0.8073    & 3.4124 & 1.1265\\
         $c_8$  & 1.7413 & 1.1265    & 3.7654 & 2.0578\\
         $c_9$  & 2.8614 & 2.0578    & 4.0119 & 2.7109\\
         $c_{10}$ & 3.5410 & 2.7109  & 4.3412 & 3.1543\\
         $c_{11}$ & 4.0034 & 3.1543  &   -     & 3.5418\\\hline
  \end{tabular}
         \caption{Values $c_i$, $0\leq i\leq 11$ in Tests 1 and 2.}    \label{tbl:cis}
\end{center}
    \end{table}

\begin{figure}
    \centering
    \includegraphics[width=0.6\textwidth]{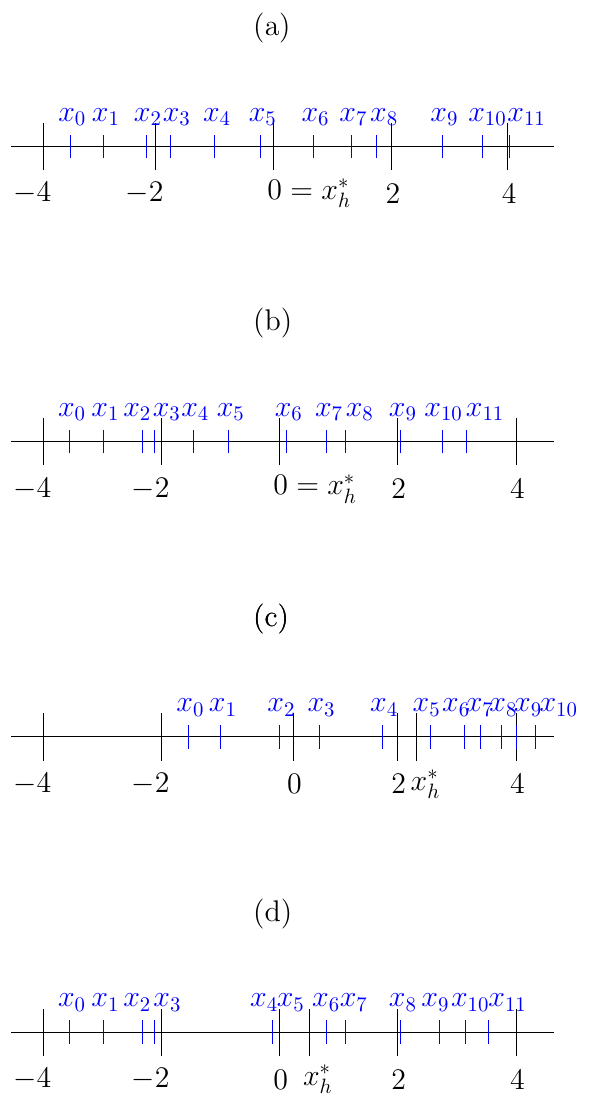}
    \caption{Representation of the different {nonuniform} meshes considered in Tests 1 and 2 for $h=1$, $z=0$ and the values of $c*$ indicated in each case. Note that to simplify the notation we have used $x_i$ to denote $x_{i,h}$. (a) Mesh for Test 1 and smooth initial data, (b) mesh for Test 1 and discontinuous initial data, (c) mesh for Test 2 and smooth initial data and (d) mesh for Test 2 and discontinuous initial data.}
    \label{fig:malles}
\end{figure}

\subsubsection{Test 1: Reconstruction of a smooth function}\label{sec:smooth_point}

Let us consider the function $f(x):=xe^x$ and the grid given in Figure \ref{fig:malles} (a), defined by the values $c_i$, $i=0,\hdots,11$ showed in Table \ref{tbl:cis} (column 1), with  $N=12$ and $z=0$. We consider the reconstruction point $x_h^*=z+c^*h$, with $c^*=0$, hence in this case $x_h^*$ is always located between the nodes $x_{5,h}$ and $x_{6,h}$.

First, we interpret the given data as point values, namely, we consider $f_{i,h}=f(x_{i,h})$, $0\leq i\leq 11$. In this case, the expected accuracy order should be the optimal twelfth-order accuracy that would attain the reconstruction procedure without WENO weights. By applying it to our WENO scheme, Table \ref{tbl:smooth_point} shows the numerical errors and accuracy order obtained, which clearly shows the convergence to the aforementioned optimal order, as we have proved in Theorem \ref{acc_standard}.

\begin{center}
    \begin{table}
        \begin{tabular}{|c|cc|cc|}       \hline
         & \multicolumn{2}{c|}{Point values} &  \multicolumn{2}{c|}{Cell averages} \\
        \hline
        $n$ &$E_n$ & $o_n$ & $E_n$ & $o_n$  \\
        \hline
        $0$  & 5.5486e-14 &         & 4.5796e-13 &      \\
        $1$  & 1.3161e-17 & 12.0416 & 2.2884e-16 & 10.9667 \\
        $2$  & 3.1728e-21 & 12.0183 & 1.1319e-19 & 10.9813 \\
        $3$  & 7.7003e-25 & 12.0085 & 5.5649e-23 & 10.9902 \\
        $4$  & 1.8746e-28 & 12.0041 & 2.7267e-26 & 10.9950 \\
        $5$  & 4.5703e-32 & 12.0020 & 1.3337e-29 & 10.9974 \\
        $6$  & 1.1150e-35 & 12.0010 & 6.5184e-33 & 10.9987 \\
        $7$  & 2.7212e-39 & 12.0005 & 3.1842e-36 & 10.9994 \\
        $8$  & 6.6426e-43 & 12.0002 & 1.5551e-39 & 10.9997 \\
        $9$  & 1.6215e-46 & 12.0001 & 7.5944e-43 & 10.9998 \\
        $10$ & 3.9587e-50 & 12.0001 & 3.7084e-46 & 10.9999\\
        $11$ & 9.6648e-54 & 12.0000 & 1.8108e-49 & 11.0000\\
        $12$ & 2.3595e-57 & 12.0000 & 8.8419e-53 & 11.0000\\
        $13$ & 5.7605e-61 & 12.0000 & 4.3173e-56 & 11.0000\\
        $14$ & 1.4063e-64 & 12.0000 & 2.1081e-59 & 11.0000\\
        $15$ & 3.4335e-68 & 12.0000 & 1.0293e-62 & 11.0000\\
        $16$ & 8.3827e-72 & 12.0000 & 5.0261e-66 & 11.0000\\
        $17$ & 2.0465e-75 & 12.0000 & 2.4541e-69 & 11.0000\\
        $18$ & 4.9964e-79 & 12.0000 & 1.1983e-72 & 11.0000\\
        $19$ & 1.2198e-82 & 12.0000 & 5.8511e-76 & 11.0000\\
        \hline
        \end{tabular}
        \caption{Test 1: Numerical errors and orders of the reconstruction of the smooth function from point values with grid given in column 1 of Table \ref{tbl:cis} and from cell average with grid given in column 2 of Table \ref{tbl:cis}.}    \label{tbl:smooth_point}
    \end{table}
\end{center}

Now, we interpret the data as cell averages, namely, we consider
$$\bar{f}_{i,h}=\frac{1}{x_{i+1}-x_i}\int_{x_i}^{x_{i+1}}f(x)dx,\quad0\leq i<11,$$
using the same function $f$ and the grid given in Figure \ref{fig:malles} (b) defined by the values $c_i$, $i=0,\hdots,11$ in Table \ref{tbl:cis} (column 2).
We consider $x_h^*=z+ch$, with $c^*=0$. As can be seen in Table \ref{tbl:smooth_point}, the behaviour of the numerical solutions is the same as before: the expected order of accuracy is reached.

\subsubsection{Test 2: Reconstruction of a discontinuous function}\label{sec:discnt_point}

Let us now consider now the function
\[f(x):=
\begin{cases}
  xe^x, & x\leq0, \\
  2xe^x+1, & x>0,
\end{cases}
\]
and the grid given by the data shown in column 3 of Table \ref{tbl:cis} with $c^*=2.3251$. Note that the reconstruction point is now located between $x_{4,h}$ and $x_{5,h}$. In this case, the discontinuity divides the stencil into two halves, where the first three nodes contain the information corresponding to the left side of the discontinuity, namely, $xe^x$, and the last eight nodes contain the information related to the right side of the discontinuity, that is, $2xe^x+1$. The expected behaviour is therefore an accuracy drop to sixth order (since $\lfloor\frac{11+1}{2}\rfloor=6$), which is confirmed in Table \ref{tbl:discnt_point}.

\begin{center}
    \begin{table}
        \begin{tabular}{|c|cc|cc|}       \hline
         & \multicolumn{2}{c|}{Point values} &  \multicolumn{2}{c|}{Cell averages} \\
        \hline
        $n$ & $E_n$ & $o_n$ & $E_n$ & $o_n$  \\
        \hline
         $0$  & 1.9479e-01 &       & 9.4953e-02 &    \\
         $1$  & 1.7007e-03 & 6.8396& 5.1777e-04 & 7.5188\\
         $2$  & 8.4877e-06 & 7.6466& 3.3859e-06 & 7.2566\\
         $3$  & 6.8988e-08 & 6.9429& 3.3852e-08 & 6.6442\\
         $4$  & 7.6213e-10 & 6.5002& 4.2252e-10 & 6.3241\\
         $5$  & 9.9623e-12 & 6.2574& 5.8992e-12 & 6.1623\\
         $6$  & 1.4219e-13 & 6.1306& 8.7129e-14 & 6.0812\\
         $7$  & 2.1227e-15 & 6.0658& 1.3235e-15 & 6.0406\\
         $8$  & 3.2417e-17 & 6.0330& 2.0391e-17 & 6.0203\\
         $9$  & 5.0075e-19 & 6.0165& 3.1638e-19 & 6.0102\\
         $10$ & 7.7796e-21 & 6.0083& 4.9261e-21 & 6.0051\\
         $11$ & 1.2120e-22 & 6.0041& 7.6836e-23 & 6.0025\\
         $12$ & 1.8911e-24 & 6.0021& 1.1995e-24 & 6.0013\\
         $13$ & 2.9528e-26 & 6.0010& 1.8734e-26 & 6.0006\\
         $14$ & 4.6121e-28 & 6.0005& 2.9265e-28 & 6.0003\\
         $15$ & 7.2051e-30 & 6.0003& 4.5722e-30 & 6.0002\\
         $16$ & 1.1257e-31 & 6.0001& 7.1437e-32 & 6.0001\\
         $17$ & 1.7588e-33 & 6.0001& 1.1161e-33 & 6.0000\\
         $18$ & 2.7481e-35 & 6.0000& 1.7440e-35 & 6.0000\\
         $19$ & 4.2939e-37 & 6.0000& 2.7249e-37 & 6.0000\\  \hline
        \end{tabular}
         \caption{Test 2: Numerical errors and orders of the reconstruction of the discontinuous function from point values with grid given in column 3 of Table \ref{tbl:cis} and from cell average with grid given in column 4 of Table \ref{tbl:cis}.}\label{tbl:discnt_point}
    \end{table}
\end{center}

We repeat the experiment for cell averages with nodes given by the values $c_i$, $0\leq i\leq 11$ in the fourth column of Table \ref{tbl:cis} and $c^*=0.5041$. In this case, the reconstruction point is located in the cell with bounds $x_{5,h}$ and $x_{6,h}$ as it can be seen in Figure \ref{fig:malles} (d). As before, the discontinuity separates the cells into two parts. The first one is given by the first five cells, and the second one contains six cells with the information related to the right side of the discontinuity. As we expected, in Table \ref{tbl:discnt_point} we observe that the accuracy order reaches sixth order, as we have studied in Theorem \ref{acc_standard}.

\subsection{Numerical experiments involving conservation laws}




To illustrate the applicability of our proposed WENO interpolator in the context of numerical solutions of PDEs we next consider some experiments using one-dimensional hyperbolic conservation laws $u_t+f(u)_x=0.$ First, we will test our method using some academic examples such as the linear advection equation, with $f(u)=u$, and Burgers equation, which corresponds with $f(u)=u^2/2$. Later, we will perform two more sophisticated experiments to compare the results obtained with the non-uniform approach and the uniform one in more challenging scenarios.

For this purpose, a finite volume scheme will be used (see e.g. \cite{shu98}), which is suitable to attain a method with high order accuracy. Given a domain $\Omega=[a,b]$ and a partition of $\Omega$, $\Big\{x_{i+\frac{1}{2}}\Big\}_{i=0}^N$, the semidiscrete {form of the conservation law reads}
\begin{equation}\label{eq:ugen}
\overline{u}_t(x_i,t)=-\frac{1}{\Delta x_i}\left(f\left(u\left(x_{i+\frac{1}{2}},t\right)\right)-f\left(u\left(x_{i-\frac{1}{2}},t\right)\right)\right),
\end{equation}
where
$$\overline{u}(x_i,t):=\frac{1}{\Delta x_i}\int_{x_{i-\frac{1}{2}}}^{x_{i+\frac{1}{2}}}u(x,t)dx,\quad x_i:=\frac{x_{i-\frac{1}{2}}+x_{i+\frac{1}{2}}}{2},\quad\Delta x_i:=x_{i+\frac{1}{2}}-x_{i-\frac{1}{2}}.$$
To approximate \eqref{eq:ugen} {for approximations $\overline{u}_{i}(t)\approx\overline{u}_t(x_i,t)$}, we use the conservative difference formula:
\begin{equation*}
{\overline{u}'_i(t)}=-\frac{1}{\Delta x_i}\left(\overline{f}_{i+\frac12}-\overline{f}_{i-\frac12}\right),
\end{equation*}
where the numerical fluxes, $\overline{f}_{i+\frac12}\approx f(u(x_{i+\frac{1}{2}},t))$,  are approximated using:

\begin{description}
\item[For the autonomous linear advection equation] Left-biased upwind reconstructions are applied:
$$\overline{f}_{i+\frac12}=q^-_{5,h}(x_{i+\frac{1}{2}}),$$
where $q^-_{5,h}$ is the WENO reconstruction using five cells left-biased, i.e., $\overline{u}_{i+l}$, $-3\leq l\leq 2$.

\item[For the rest of the numerical experiments in this work] Local {Lax-Friedrichs} flux splitting is utilized:
 $$\overline{f}_{i+\frac12}=\frac{1}{2}\left(f(q^+_{5,h}(x_{i+\frac{1}{2}}))+f(q^-_{5,h}(x_{i+\frac{1}{2}}))
-\alpha(q^+_{5,h}(x_{i+\frac{1}{2}})-q^-_{5,h}(x_{i+\frac{1}{2}}))\right),$$
where $q^{\pm}_{5,h}$ are the WENO reconstructions using five cells left and right-biased and $\alpha=\max_{u} |f'(u)|$. 
\end{description}

The time-stepping solver that will be used along all the experiments is the second- or third-order Runge-Kutta TVD scheme \cite{ShuOsher89}. The spatial discretization will always have order 5.

\subsubsection{Construction of the {nonuniform} grid}\label{sec:nugrid}
For the sake of reproducibility, we provide all the details about the non-uniform grid construction in each experiment. Given $n\in\N$, we define
$$x_j:=-1+\frac{2j}{n}+\frac{2}{n}R(j, \xi),\quad 0\leq j\leq n,$$
with $R(0,\xi)=R(n,\xi)=0$ and
$$R(j,\xi):=-\xi-2\xi r(s_1(j-1),s_2(j-1),s_3(j-1)),$$
for $1\leq j\leq n-1$. The parameter $\xi$ controls the fluctuation of the cell interfaces with respect to a uniform grid, and the function $r:\N^3\to [0,1]$ is the Wichmann-Hill pseudo-random number generator given by
\begin{align*}
    r(s_1(j-1),s_2(j-1),s_3(j-1))=\frac{s_1(j)}{30269}+\frac{s_2(j)}{30307}+\frac{s_3(j)}{30323}\mod 1,
\end{align*}
with seed update
\begin{align*}
    s_1(j)&=171s_1(j-1)\mod 30269,\\
    s_2(j)&=172s_2(j-1)\mod 30307,\\
    s_3(j)&=170s_3(j-1)\mod 30323.
\end{align*}
In all the experiments, the selected starting seeds will be chosen as
\begin{align*}
    s_1(0)&=874,\\
    s_2(0)&=1421,\\
    s_3(0)&=957.
\end{align*}
In the experiments analyzing the order under smooth solutions, namely, those involving several simulations in which the number of cells is duplicated, the seeds are not restarted upon the next doubled resolution; instead, the last updated seed from the previous experiment with halved resolution is used as starting seed for the next one.

\subsubsection{Test 3: Linear advection equation}\label{sec:linadv_smooth}

First of all, we compute the approximations of the solution of the linear advection equation:
$$u_t+u_x=0, (x,t)\in[-1,1]\times[0,\infty),$$
with initial and boundary conditions:
\begin{equation*}
\begin{split}
&u(x,0)=0.25+0.5\sin(\pi x), \,\, x\in [-1,1],\\
&{u(-1+x,t)=u(1+x,t)}, \,\, t\in[0,\infty),
\end{split}
\end{equation*}
which exact solution is given by:
$$u(x,t)=0.25+0.5\sin\left(\pi(x-t)\right).$$

As the method used to solve the ODE is third-order accurate, {for test purposes} we choose $\Delta t=(\widetilde{\Delta x})^{5/3}$ with
$\widetilde{\Delta x}=\min_{i}\{\Delta x_i\}$. With this selection, we guarantee that $\Delta t^3=O(\Delta x^{5})$, with ${\Delta x}=\max_{i}\{\Delta x_i\}$ and $\Delta t/\Delta x \leq 1$.
In this case, the fluctuation parameter is $\xi=0.1$. We calculate the error for $T=1$ and the numerical order in a similar way as we computed in Equation \eqref{order} for values of $n=20\cdot 2^j$ with $j=0,\hdots,4$. The results are displayed in Table \ref{tbl:linadv}. It is clear that the expected order of accuracy is reached.

\begin{center}
    \begin{table}
        \begin{tabular}{|c|cc|cc|}
        \hline
        $n$ & Error $\|\cdot\|_1$ & Order $\|\cdot\|_1$ & Error $\|\cdot\|_{\infty}$ & Order $\|\cdot\|_{\infty}$ \\
        \hline
        $20$ & 5.55e-05 & - & 8.84e-05 & - \\
        $40$ & 1.79e-06 & 4.96 & 9.32e-05 & 4.95 \\
        $80$ & 5.63e-08 & 4.99 & 9.32e-08 & 4.94 \\
        $160$ & 1.77e-09 & 4.99 & 2.94e-09 & 4.99 \\
        $320$ & 5.57e-11 & 4.99 & 9.31e-11 & 4.98 \\
        \hline
        \end{tabular}
        \caption{Test 3: Errors and orders obtained for linear advection equation with smooth initial data with $n=20\cdot 2^j$, $j=0,\hdots,4$}\label{tbl:linadv}
    \end{table}
\end{center}

\begin{figure}
    \centering
    \begin{tabular}{cc}
    (a) & (b)\\
    \includegraphics[width=0.5\textwidth]{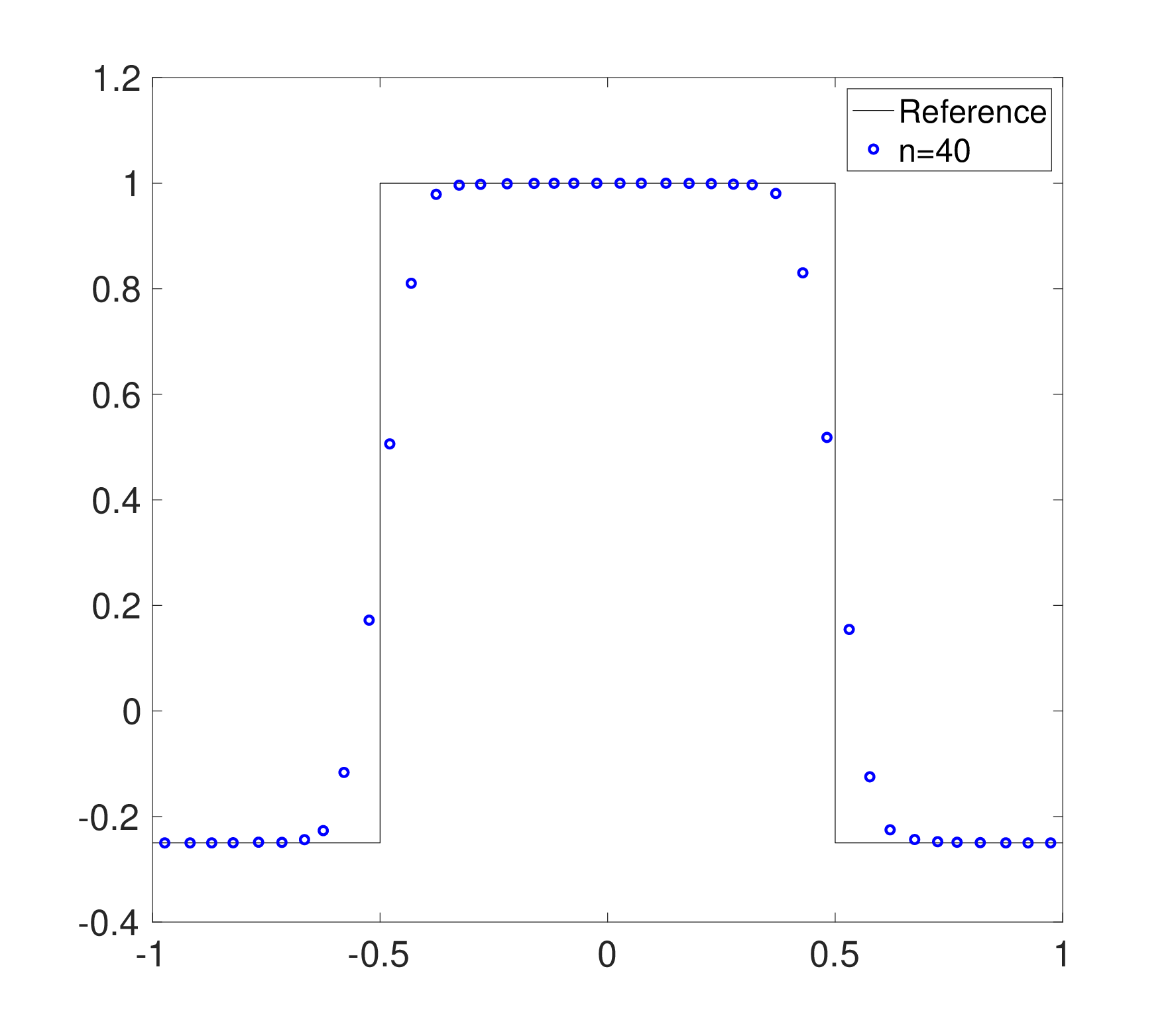} &
\includegraphics[width=0.5\textwidth]{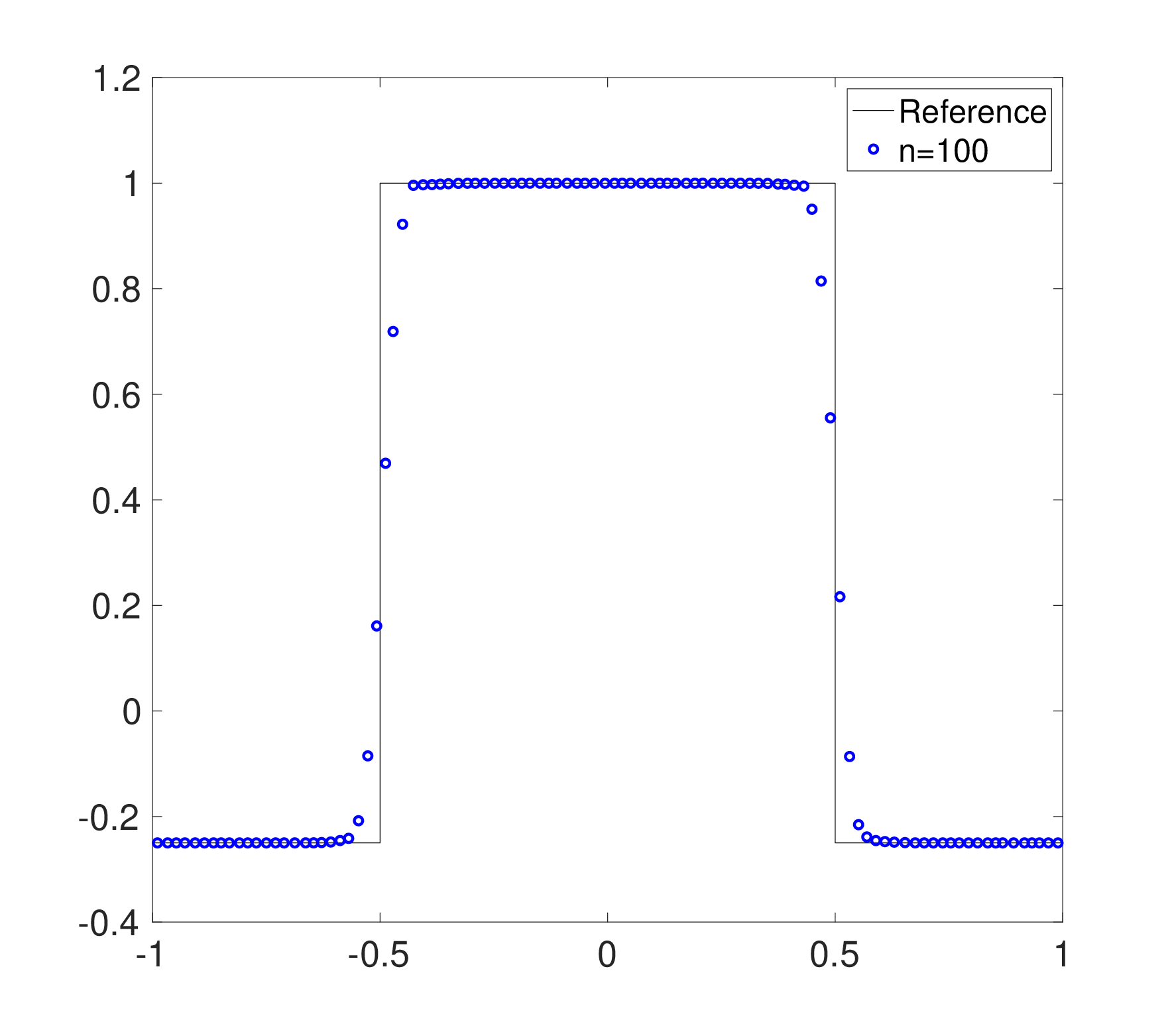}
\end{tabular}
    \caption{Test 3: Numerical solution of the linear advection equation with discontinuous initial data for (a) $n=40$, (b) $n=100$}
    \label{fig:linadv}
\end{figure}

We consider now the discontinuous initial data given by:
\[u(x,0)=
\begin{cases}
    -0.25, & x\leq0, \\
    1, & x>0,
\end{cases}
\]
$\Delta t=0.9\widetilde{\Delta x}$ and
the fluctuation parameter $\xi=0.25$. The simulation is now run until $T=1.5$ and with $n=40$ and $n=100$. We can see in Figure \ref{fig:linadv} that the behaviour of the numerical solution is oscillation-free and does not present Gibbs phenomenon close to the discontinuities present in the solution.

\subsubsection{Test 4: Burgers equation}

\begin{figure}
    \centering
    \begin{tabular}{cc}
    (a) & (b)\\
    \includegraphics[width=0.5\textwidth]{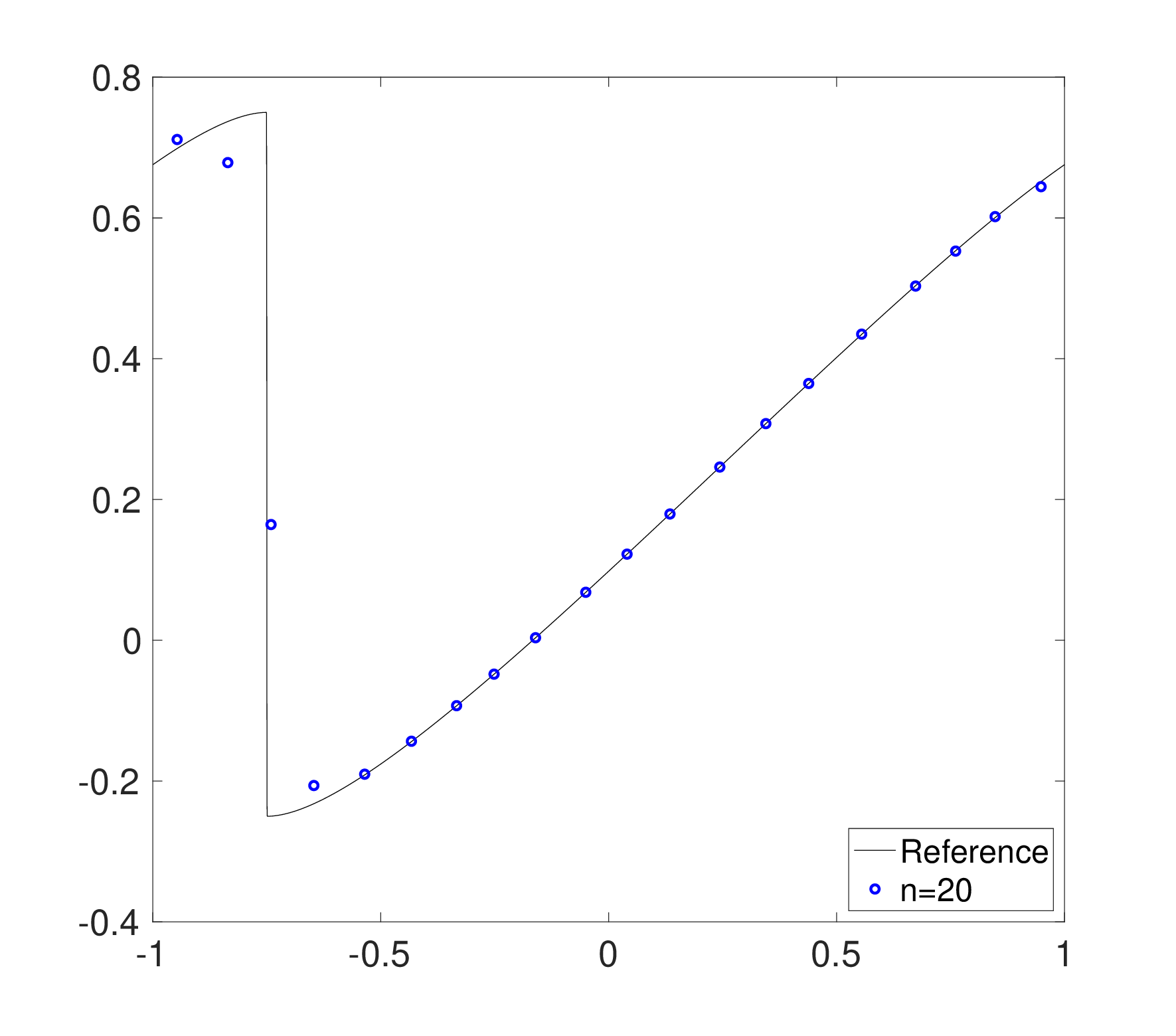} &
\includegraphics[width=0.5\textwidth]{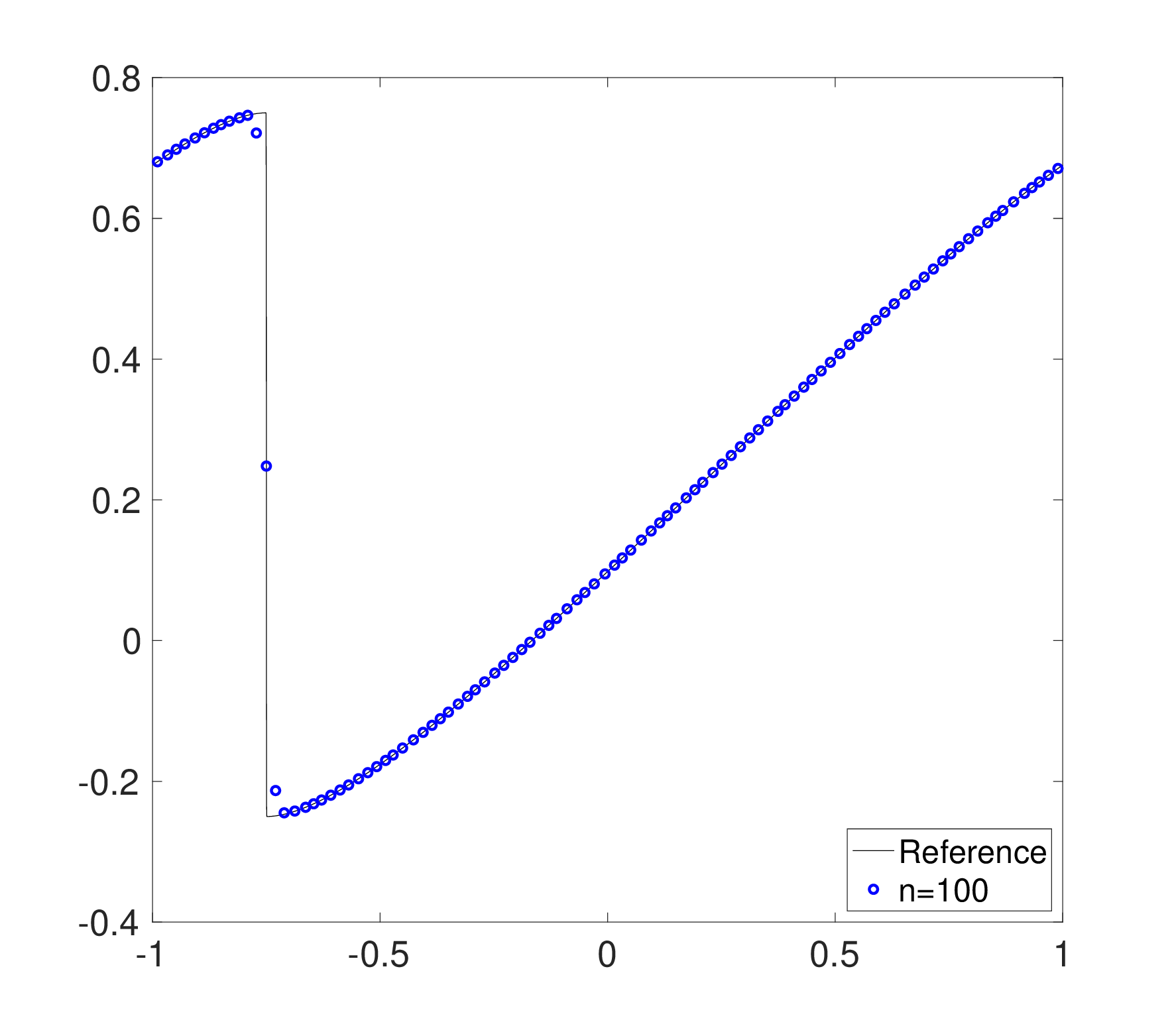}
\end{tabular}
    \caption{Test 4: Numerical solution of the Burgers equation with discontinuous initial data: (a) $n=20$, (b) for $n=100$}
    \label{fig:burgers}
\end{figure}

If we repeat the same experiments as in Test 3 but for the Burgers equation with $T=0.3$ for smooth initial data and $T=1$ for non-smooth initial data, we obtain
the results shown in Table \ref{tbl:burgers} and in Figure \ref{fig:burgers}. Again, the expected approximation order, 5, is achieved for both norm 1 and norm infinity. The simulations in Figure \ref{fig:burgers} show that for small values of $n$ some numerical oscillations appear near the boundaries. However, this problem is overcome when we refine the mesh.


\begin{center}
    \begin{table}
        \begin{tabular}{|c|cc|cc|}
        \hline
        $n$ & Error $\|\cdot\|_1$ & Order $\|\cdot\|_1$ & Error $\|\cdot\|_{\infty}$ & Order $\|\cdot\|_{\infty}$ \\
        \hline
        $40$ & 1.21e-05 & - & 1.57e-04 & - \\
        $80$ & 4.18e-07 & 4.86 & 5.73e-06 & 4.77 \\
        $160$ & 1.16e-08 & 5.17 & 1.63e-07 & 5.14 \\
        $320$ & 3.45e-10 & 5.07 & 5.46e-09 & 4.90 \\
        $640$ & 1.09e-11 & 4.99 & 2.06e-10 & 4.73 \\
        \hline
        \end{tabular}      
        \caption{Test 4: Errors and orders obtained for Burgers equation with smooth initial data with $n=20\cdot 2^j$, $j=0,\hdots,4$}        \label{tbl:burgers}
    \end{table}
\end{center}

%
%
%
%

\subsubsection{Test 5: Scalar equation with solution approaching a Dirac delta function}\label{sec:delta}

It is the goal of this test to show a numerical experiment where the results obtained using the nonuniform approach significantly outperform those of the uniform one. We consider the equation (see \cite{MuletVecil})
\begin{equation}\label{eq:delta}
    u_t+\left(\sin(x) u\right)_x=0,\quad x\in[0,2\pi],
\end{equation}
with initial condition $u(x,0)=1$, $x\in[0,2\pi]$, and periodic boundary conditions. By the method of characteristics, the exact solution of \eqref{eq:delta} is given by
$$u(x,t)=\frac{2e^t}{1-\cos(x)+\left(1+\cos(x)\right)e^{2t}}.$$ 
For a fixed value of $x\in[0,2\pi]$, there holds
$$\lim_{t\to\infty}u(x,t)=\begin{cases}
    +\infty & x=\pi,\\
    0 & x\neq\pi,
  \end{cases}$$
  and for fixed $t\geq 0$
  $$
  \int_{0}^{2\pi} u(x, t)dx=2\pi,
  $$
  therefore, the solution of this equation behaves progressively like a Dirac delta function as time increases. This means that the solution will have a strong gradient near $x=\pi$, while it would remain close to zero in the rest of the domain. The structure of the solution suggests that the use of a nonuniform grid with decreasing distance between nodes as we approach the point $x=\pi$ could be convenient.

We set $T=8$ and consider uniform grids with $n=100 \cdot 2^{j}-1, j=1,\dots,10$ cells. We define non-uniform grids of $2m$ cells with cell boundaries defined by
\begin{align*}
x_{\mig}&=0,\\ 
x_{j+\mig}-x_{j-\mig}&=\kappa(x_{j+\mig+1}-x_{j+\mig}),\quad j=1,\dots,m-1,\\ 
x_{m+\mig}&=\pi,
\end{align*}
and, symmetrically, $x_{m+j+\mig}=2\pi-x_{m-j+\mig}$, $j=1,\dots,m$, with $\kappa>1$.

 The errors are computed as
  \begin{align}\label{error:delta}
    \sum_{j=1}^{n} (x_{j+\mig}-x_{j-\mig})
    \left|\frac{1}{x_{j+\mig}-x_{j-\mig}}\int_{x_{j-\mig}}^{x_{j+\mig}}u(x, T)dx - u_{j}^{N}\right|,
  \end{align}
  where $u_{j}^{N}$ is the result of $N$ time steps of length $\Delta t=T/N$, subject to the following CFL restriction
  \begin{align*}
  \frac{\Delta t}{\min\limits_{j} x_{j+\mig}-x_{j-\mig}} \approx 0.8.
  \end{align*}
  The ODE solver used in this test is TVD-RK2.

\begin{figure}
    \centering
    \begin{tabular}{cc}
    (a) & (b)\\
    \includegraphics[width=0.5\textwidth]{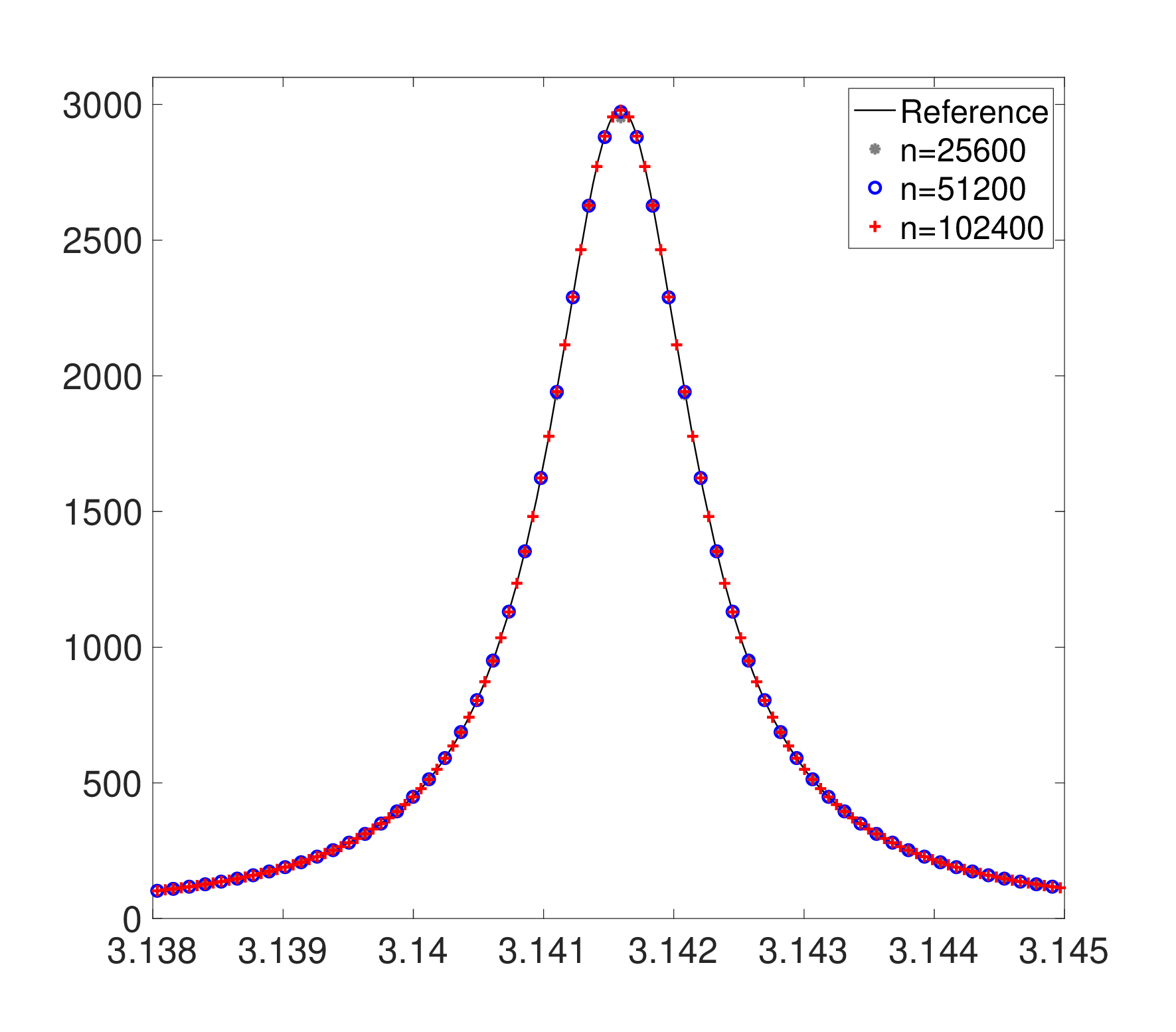} &
    \includegraphics[width=0.5\textwidth]{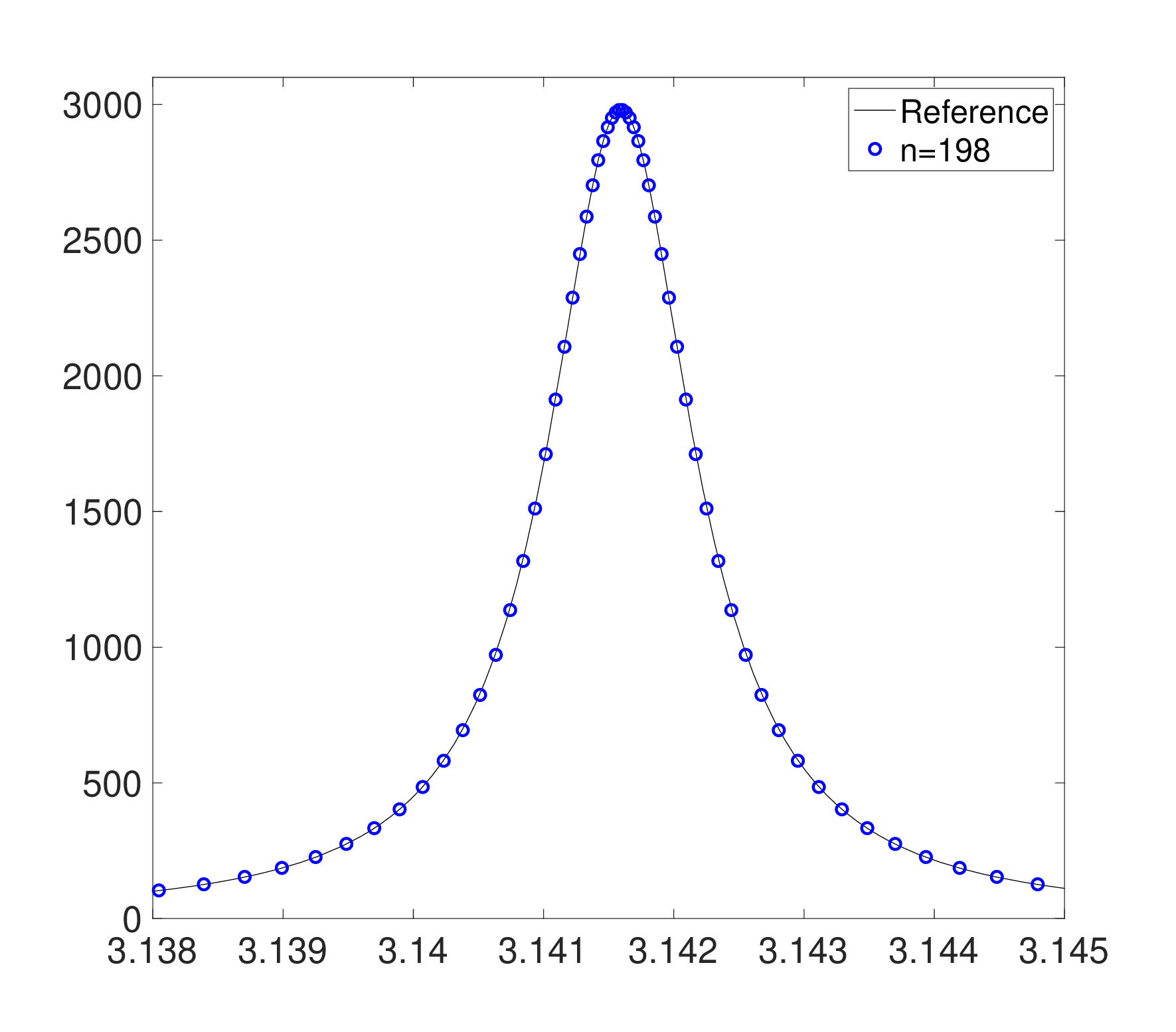}\\
    (c) & (d)\\
    \includegraphics[width=0.5\textwidth]{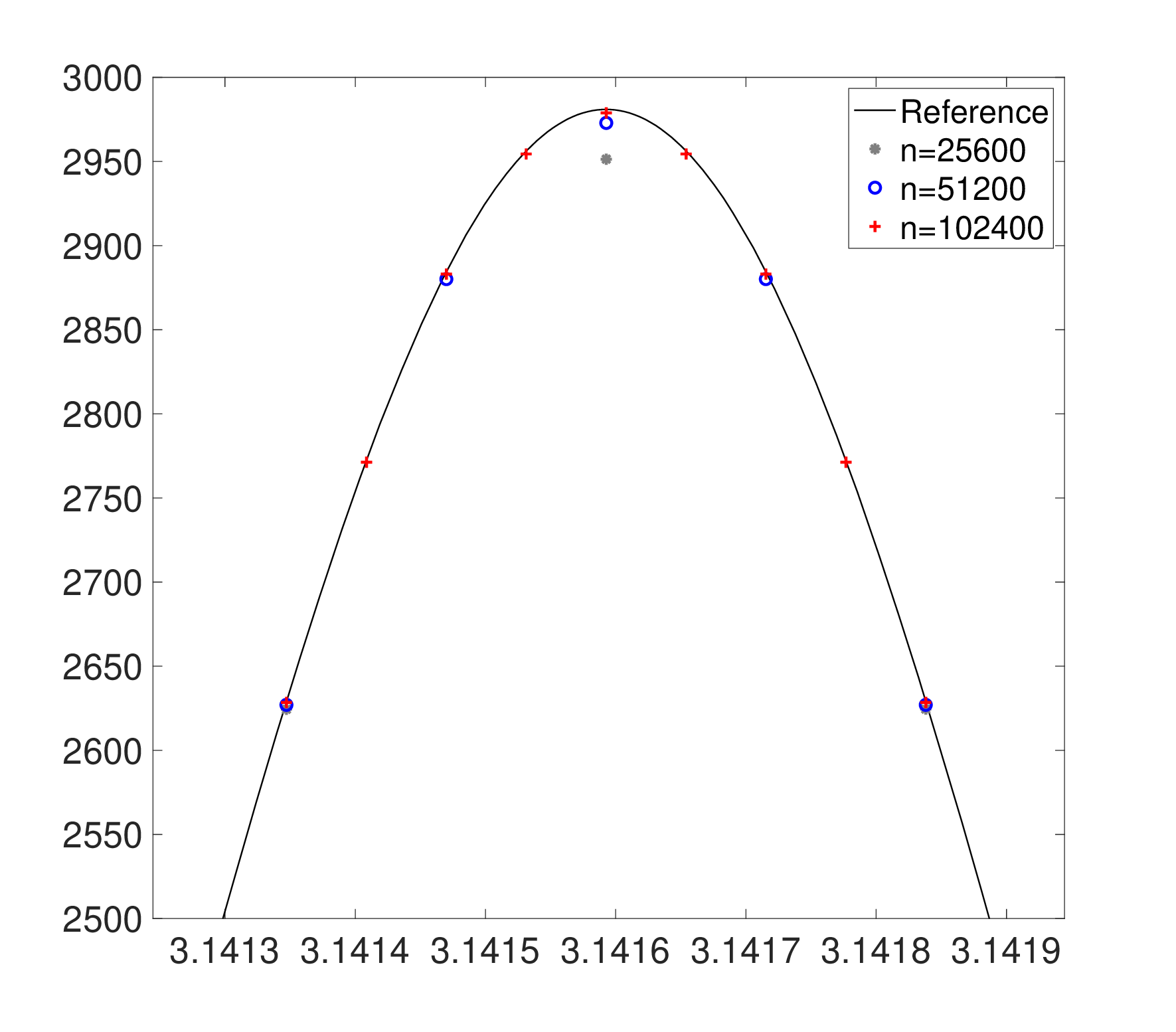} &
    \includegraphics[width=0.5\textwidth]{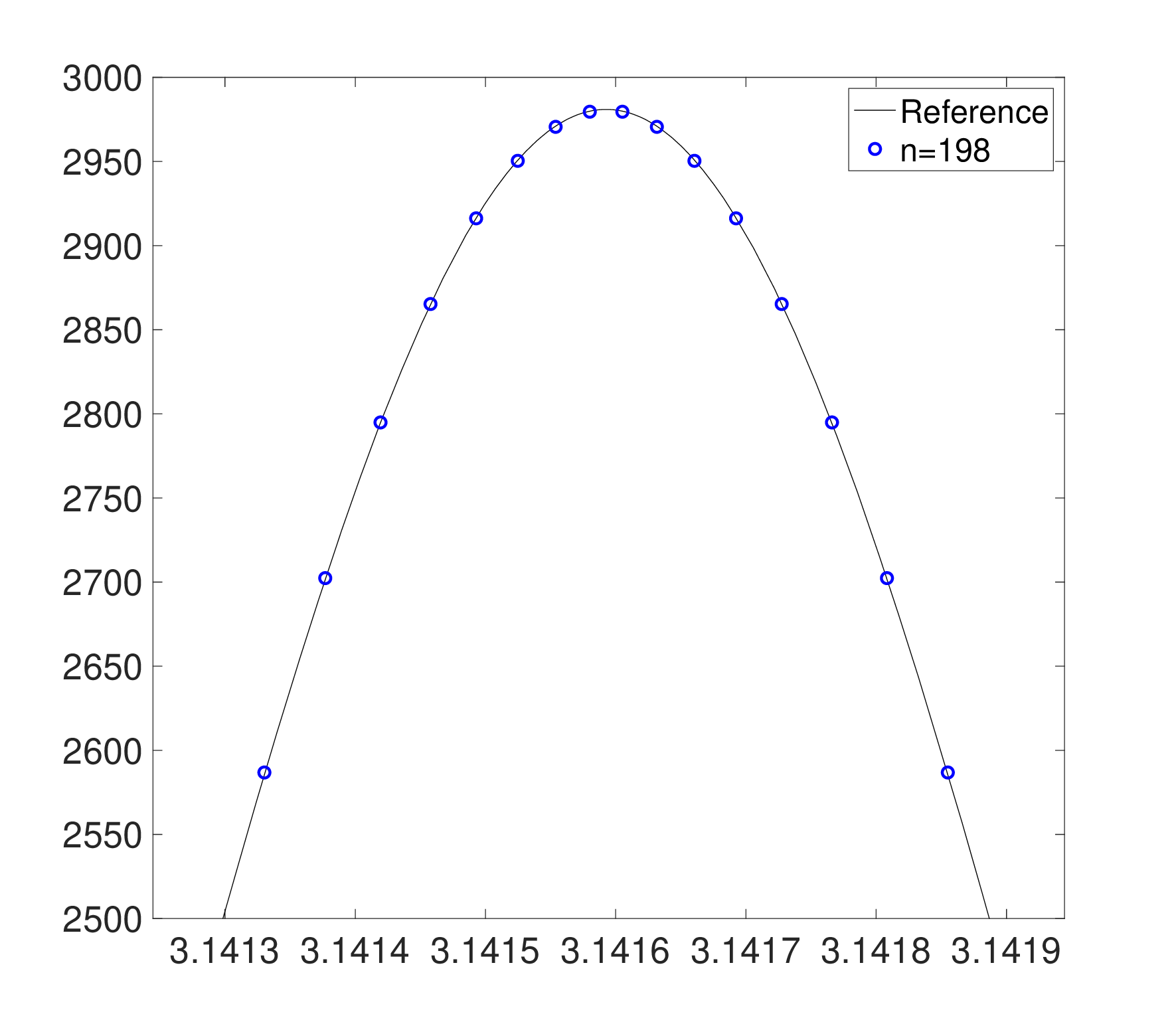} 
\end{tabular}
    \caption{Test 5: Numerical solution of the scalar equation converging to a Dirac delta function in the interval $[3.138,3.145]$, obtained using: (a) uniform meshes with $n=25599, 51199, 102399$, (b) non-uniform mesh with $n=198$ and $q=1.1$. (c) and (d) are enlarged views of (a) and (b) respectively.}
    \label{fig:delta}
  \end{figure}

In Figure \ref{fig:delta} we show the numerical solutions obtained using both uniform and non-uniform meshes. Figures ~\ref{fig:delta} (c) and (d) show that the accuracy of the approximated solutions obtained using non-uniform grids is higher than the ones obtained using uniform grids, even when a significantly lower resolution is considered.
  
Table \ref{tbl:delta} displays the errors and orders obtained when using uniform grids. For comparison, the error obtained using a non-uniform grid with $2m=198$ cells and $q=1.1$ is $2.87e-2$, which is comparable to the error obtained with a uniform grid of 25599 cells. If we consider $2m=398$ cells and $q=1.04$ the computed error is $2.67e-3$, which in this case is comparable to the error obtained with a uniform grid with 51199 cells.

It is worth mentioning that although the number of cells on the non-uniform grids considered is much smaller than the number of cells on the uniform ones, the smallest cell size for a non-uniform grid with $2m=198$ cells and $q=1.1$ is $2.51e-5$ while for the grid with $2m=398$ and $q=1.04$ is $5.13e-5$. Those sizes are smaller than the cell size for the uniform grids with 25599 and 51199 cells, given by $2.45e-4$ and $1.23e-4$, respectively.

\begin{table}
  \begin{tabular}{|c|c|c|}
    \hline
      $n$& Error      & Order\\\hline
     199 & 2.03e+00 & - \\ 
     399 & 2.09e+00 & 9.7043e-01\\ 
     799 & 2.08e+00 & 1.0049e+00\\ 
    1599 & 1.84e+00 & 1.1325e+00\\ 
    3199 & 1.17e+00 & 1.5657e+00\\ 
    6399 & 3.47e-01 & 3.3787e+00\\ 
   12799 & 5.92e-02 & 5.8592e+00\\ 
   25599 & 2.02e-02 & 2.9247e+00\\ 
   51199 & 3.17e-03 & 6.3790e+00\\ 
  102399 & 7.12e-04 & 4.4556e+00\\ \hline
\end{tabular}
\caption{Test 5: Errors and orders obtained using uniform grids with $n$ cells, $n=100\cdot 2^j -1$, $j=1,\hdots,10$. The errors are computed using formula \eqref{error:delta}.}
\label{tbl:delta}
\end{table}

\subsubsection{Test 6: Shu-Osher problem}\label{sec:shuosher}
In this test, we compare the performance of our method against the CWENO-type methods. Given that the method proposed by Levy, Puppo and Russo \cite{LevyPuppoRusso} can have negative weights depending on the distribution of the grid and the relative position of the reconstruction point, we will use the alternative proposed in \cite{BaezaBurgerMuletZorio2019}, which was designed to avoid this issue as one of its purposes. The main goal of this test is to compare the impact of the computation of the smoothness indicators in the global cost of the method, which has been discussed theoretically in Subsection \ref{sec:comp}.

We start by describing the model of the 1D Euler equations for gas dynamics. If we denote the density as $\rho$,  the velocity as $v$ and the specific energy of the system by $E$, then the equations are given by $\boldsymbol{u} = ( \rho, \rho v, E)^{\mathrm{T}}$ and $\boldsymbol{f} (\boldsymbol{u}) =  (\rho v, p+\rho v^2, v(E+p))^{\mathrm{T}}$. The variable~$p$ represents the pressure and is defined by the equation of state
$$p=\left(\gamma-1\right)\left(E-\frac{1}{2}\rho v^2\right),$$
where $\gamma$~is the adiabatic constant taken as $\gamma =1.4$. We consider as spatial domain $\Omega=(-5,5)$, and as initial condition
\begin{align*} 
(\rho,v,p) (x, 0) = 
\begin{cases} \displaystyle 
  \biggl(\frac{27}{7}, 
    \frac{4\sqrt{35}}{9}, 
    \frac{31}{3} \biggr) & \text{if  $x\leq-4$,}  \\[4mm]  \displaystyle
   \biggl(1+\frac{1}{5}\sin(5x), 0, 1 \biggr) 
    & \text{if $x>-4$.} 
\end{cases} \end{align*}
This initial condition specifies the interaction of a Mach~3 shock with a sine wave and is accompanied by left inflow and right outflow boundary conditions.

\begin{figure}
\centering
\includegraphics[width=\textwidth]{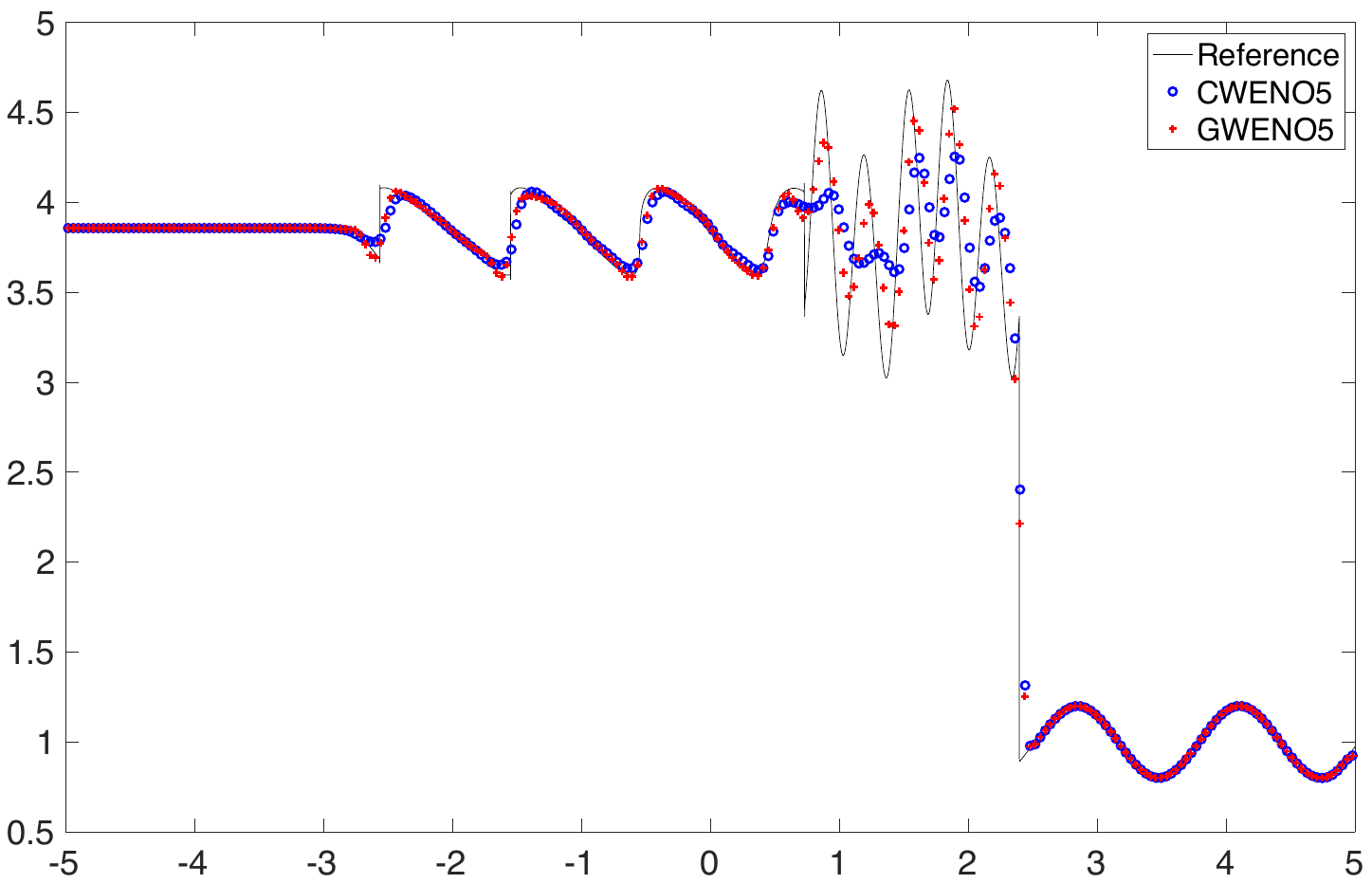}
\caption{Test 6: Numerical solution of the Shu-Osher problem with CWENO-type schemes with uniform (blue) and non-uniform (red) grids with $n=256$ cells.}
\label{fig:shuosher}
\end{figure}

 
We simulate until $T=1.8$ with a resolution of $n=256$ cells, $\text{CFL}=0.5$ and the same parameters as those stated in Section \ref{sec:nugrid} to generate the non-uniform grid. Then, we compare the results against a reference solution computed with a resolution of $N=65536$ cells. The results of the experiment can be seen in Figure \ref{fig:shuosher}.

In this simulation, our method is 1.13 times faster than CWENO for the fifth-order scheme and 8.24 times faster for the seventh-order scheme, mostly due to the considerably more complex computations required to obtain the classical Jiang-Shu smoothness indicators on a nonuniform grid, required by the chosen CWENO method. Furthermore, our method produces a better approximation than the CWENO method, as Figure \ref{fig:shuosher} reflects.

\section{Conclusions}\label{conclusions}

In this paper, we have presented a novel WENO approach capable of performing essentially non-oscillatory reconstructions in a general context of {nonuniform} grids. Both the theoretical results and the numerical experiments show that the scheme has the desired accuracy.

The efficiency of the algorithm relies strongly on the simplification of the smoothness indicators that can be performed in weight designs akin to the one proposed by Yamaleev-Carpenter \cite{YamaleevCarpenter2009}, as shown in \cite{BBMZ18}. {Another} important factor is that the coefficients involving the grid spacing can be computed and stored previously if the numerical simulation uses a fixed grid, yielding in that case a similar efficiency to the uniform grid case. As for the cases in which the grid is not fixed, then Newton polynomials, together with the associated computation of the divided differences properly scaled, are recommended instead.

Our next goal is to tackle the challenge of redesigning the weights of the method so that it can be unconditionally optimal, regardless of the order of the critical point to which the stencil converges, for stencils with four or more nodes, similarly as done in \cite{BaezaBurgerMuletZorio2018} and \cite{BaezaBurgerMuletZorioO32018} for the case of uniform stencils.

Another of our future research projects involves a generalization of the current approach into an interpolator capable of preserving the essentially non-oscillatory properties regardless of the relative position of the point to be interpolated or extrapolated with respect to the stencil. This approach can be especially interesting to handle boundary conditions in finite-difference schemes, among other possible applications.

\section*{Declarations}

\subsection*{Conflict of interest} The authors declare that they have no conflict of interest.

\subsection*{Data Availability Statements} Data sharing not applicable to this article as no datasets were generated or analysed during the current study.

\section*{{Acknowledgments}}

All authors are supported by Spanish MINECO project PID2020-117211GB-I00 and MCM, PM and DFY are also supported by GVA project CIAICO/2021/227.


\end{document}